\documentclass[a4paper]{article}

\usepackage[latin1]{inputenc}
\usepackage[T1]{fontenc}
\usepackage[francais,english]{babel}
\usepackage{tikz}
 \usetikzlibrary{arrows,shapes,trees,snakes}
\usepackage{amsmath,amssymb,amsfonts,amsthm}
\usepackage{amsfonts}
\usepackage{amssymb}
\usepackage{latexsym}
\usepackage{comment}

\def\et0{e^{tA}x_0}

\newtheorem{thm}{Theorem}
\newtheorem{Prop}{Proposition}
\newtheorem{Def}{Definition}

\newtheorem{rk}{Remark}

\author{
K. \textsc{Beauchard} 
\footnote{CMLS, Ecole Polytechnique, 91128 Palaiseau Cedex, France,
email: Karine.Beauchard@math.polytechnique.fr (corresponding author)},
P. \textsc{Cannarsa}
\footnote{Universit\`a di Roma Tor Vergata, via della Ricerca Scientifica 1, 00133, Roma, Italy, 
email: cannarsa@axp.mat.uniroma2.it},
M. \textsc{Yamamoto}
\footnote{Department of Mathematical Sciencies, The University of Tokyo, Komaba Meguro Tokyo 153-8914 (Japan)
email: myama@ms.u-tokyo.ac.jp}
\thanks{This research was performed in the framework of the GDRE CONEDP, and was partially supported by the Istituto Nazionale di Alta Matematica ``F. Severi'', \`Ecole Polytechnique, and University of Tokyo.  The authors wish to express their gratutude to the above institutions. The first author was also supported by 
the ``Agence Nationale de la Recherche'' (ANR) Projet Blanc EMAQS number ANR-2011-BS01-017-01.}
}

\title{Inverse source problem and null controllability for multidimensional  parabolic operators\\ of Grushin type}
\date{}

\begin{document}

\maketitle

\begin{abstract}
The approach to Lipschitz stability for uniformly parabolic equations introduced  by Imanuvilov and Yamamoto in 1998, based on Carleman estimates, seems hard to apply to the case of Grushin-type operators of interest to this paper. Indeed, such estimates  are still missing for parabolic operators degenerating in the interior of the space domain.  Nevertheless, we are able to prove Lipschitz stability results for inverse source problems for such operators, with locally distributed measurements in arbitrary space dimension. For this purpose, we follow a mixed strategy which combines the appraoch due to Lebeau and Robbiano, relying on Fourier decomposition, with  Carleman inequalities for 
heat equations with nonsmooth coefficients (solved by the Fourier modes).
As a corollary, we obtain a direct proof of the observability of multidimensional Grushin-type parabolic equations, 
with locally distributed observations---which is equivalent to null controllability with locally distributed controls.
\end{abstract}




\section{Introduction}

\subsection{Main results}
The relevance of subriemannian structures and hypoelliptic operators to quantum mechanics has long been acknowledged---at least since the work of Weyl~\cite{Weyl} in 1931. More recently, deep connections have been pointed out between the properties of subriemannian operators, such as the Heisenberg laplacian, and other topics of interest to current  mathematical research such as control theory and isoperimetric problems (see, for instance, \cite{capogna}).

In this paper, rather than the Heisenberg laplacian, we will study a simpler example of such type of equations, which we believe could serve as a model problem for a more general theory, in addition to being of interest in its own right. More precisely, 
we will  be concerned with equations of the form
\begin{equation} \label{Grushin}
\left\lbrace \begin{array}{ll}
\partial_t u - \Delta_x u - |x|^{2\gamma} b(x)  \Delta_y u = g(t,x,y) & (t,x,y) \in (0, \infty)\times\Omega\,,\\
u(t,x,y)=0 & (t,x,y) \in (0, \infty)\times\partial \Omega\,,
\end{array}\right.
\end{equation}
where  
$\Omega:=\Omega_1 \times \Omega_2$, 
$\Omega_1$ is a bounded open subset of $\mathbb{R}^{N_1}$, with $C^4$ boundary, such that $0 \in \Omega_1$,
$\Omega_2$ is a bounded open subset of $\mathbb{R}^{N_2}$, with $C^2$ boundary,
$N_1, N_2 \in \mathbb{N}^*:= \{1,2,3, ....\}$,
$b \in C^1(\overline{\Omega_1};(0,\infty))$, 
$\gamma \in (0,1]$  and $|.|$ is the Euclidean norm on $\mathbb{R}^{N_1}$.
In this article $\gamma \in (0,1]$ may change, but the function $b$ is fixed.

We are interested in the following questions.
\begin{itemize}
\item \underline{The inverse source problem}: for $g(t,x,y)$ given by $R(t,x)f(x,y)$, 
is it possible to recover the source term $f$ from a measurement of  
$\partial_t u |_{(T_1,T_2)\times \omega}$, 
where $\omega$ is a nonempty open subset of $\Omega$ and $R$ is suitably given?
\item \underline{The null controllability problem}:
is it possible to stear the solution to zero by applying an appropriate control
$g(t,x,y)=u(t,x,v) 1_\omega(x,v)$ localized on an open subset $\omega$ of $\Omega$?
\end{itemize}

To the best of our knowledge, there are no results on inverse 
source problems for Grushin-type equations. As for null 
controllability, the current literature  just  seems  to concern the two-dimensional case (see \cite{Grushin}).

First, we recall well-posedness and regularity results for such equations.
To this aim, we introduce the space $H^1_\gamma(\Omega)$,
which is the closure of $C^\infty_0(\Omega)$ for the topology defined by the norm
$$\|f\|_{H^1_\gamma}:= \left( \int_{\Omega} \left( |\nabla_x f|^2 + |x|^{2\gamma} |\nabla_y f|^2 \right) dx dy \right)^{1/2},$$
and the Grushin operator $G_\gamma$ defined by
$$\begin{array}{ll}
&D(G_\gamma):=  
\{ f \in H^1_\gamma(\Omega); \exists c>0 \text{ such that }\\
&\left\vert 
\int_{\Omega} \left( \nabla_x f \cdot\nabla_x g 
+ |x|^{2\gamma} \nabla_y f\cdot
\nabla_y g \right) dx dy \right\vert \leqslant c \|g\|_{L^2(\Omega)}
\quad \text{for all $g \in H^1_\gamma (\Omega)$}\},
\end{array}
$$
$$
G_\gamma u := - \Delta_x u - |x|^{2\gamma} b(x) \Delta_y u.
$$

\begin{Prop} \label{Prop:WP}
For every $u_0 \in L^2(\Omega)$ and $g \in L^2((0,T)\times \Omega)$,
there exists a unique weak solution $u \in C^0([0,T];L^2(\Omega)) \cap L^2(0,T;H^1_\gamma(\Omega))$ of (\ref{Grushin}) such that
\begin{equation} \label{IC}
\begin{array}{ll}
u(0,x,y)=u_0(x,y) & (x,y) \in \Omega.
\end{array}
\end{equation}
Moreover, $u \in C^0((0,T];D(G_\gamma))$.
\end{Prop}

We refer to \cite{Grushin} for the proof with $N_1=N_2=1$; the general case 
can be treated similarly.

\subsubsection{Inverse source problem}

Taking a source term of the form
\begin{equation} \label{form_g}
g(t,x,y)=R(t,x)f(x,y) \text{ where } R \in C^0([0,T]\times\Omega_1) \text{ and } f \in L^2(\Omega)
\end{equation}
and we will obtain Lipschitz stability estimates  for \eqref{Grushin} in the following sense.

\begin{Def}[Lipschitz stability]
Let $T>0$, $0\leqslant T_0 < T_1 \leqslant T$ and let $\omega$ be an open subset of $\Omega$.
We say that system (\ref{Grushin}) satisfies a \emph{Lipschitz stability estimate on $(T_0,T_1)\times\omega$} 
if there exists $C>0$ such that, for every $f \in L^2(\Omega)$ and $u_0 \in L^2(\Omega)$
the solution of (\ref{Grushin})(\ref{IC}) satisfies
$$\int\limits_{\Omega} |f(x,y)|^2 dx dy \leqslant C \left(
\int\limits_{T_0}^{T_1} \int\limits_{\omega} |\partial_t u(t,x,y)|^2 dx dy dt +
\int\limits_{\Omega} |G_\gamma u(T_1,x,y)|^2 dx dy
\right).$$
\end{Def}

As is easily seen from the right-hand side of the above inequality,
for the inverse problem we treat in this paper, measurements 
in $\omega$ is are taken before time $T_1$, and spatial data are measured
 in $\Omega$ at the same time $T_1$.
Moreover, we  require
the known factor $R(t,x)$ to be independent of $y$. Our method, however, could  be easily adapted to recover
a source $f=f(x)$ for $R=R(t,x,y)$.

When $\omega$ is a strip, parallel to the $y$-axis, we obtain Lipschitz stability under general assumptions on $R$.

\begin{thm} \label{Main-thm-strip}
Assume $\omega=\omega_1 \times \Omega_2$ where $\omega_1$ is an open subset of $\Omega_1$.
Suppose further that 
\begin{equation} \label{Hyp:R}
\begin{array}{c}
R, \partial_t R \in C^0([0,T]\times\overline{\Omega_1}) \text{ and } \\
\text{there exist $T_1 \in (0,T]$ and $R_0>0$ such that $R(T_1,x) 
\geqslant R_0, \forall x \in \Omega_1$}.
\end{array}
\end{equation}
\begin{enumerate}
\item If $\gamma \in (0,1)$, then system (\ref{Grushin}) satisfies the Lipschitz stability estimate 
on $(T_0,T_1)\times\omega$ for every $T_0 \in [0, T_1)$.
\item If $\gamma=1$, then there exists $T^*>0$ such that system (\ref{Grushin}) 
satisfies the Lipschitz stability estimate on $(T_0,T_1)\times\omega$ for every $T_0 \in [0,T_1 -T^*)$.
\end{enumerate}
\end{thm}

\begin{rk}
In Theorem~\ref{Main-thm-strip} above,  $T_1-T_0$ is assumed to be sufficiently large when $\gamma=1$. Indeed, in this case, the validity of a Lipschitz stability estimate on $(T_0,T_1)\times\omega$ is an open problem for a general  $T_0\in[0,T_1)$ even when $\omega$ is a strip, that is, $\omega=\omega_1 \times \Omega_2$.  On the other hand, it is known that Grushin's operator fails to be observable in arbitrary time, as an example from \cite{Grushin} shows. However, such a counterexample does not apply to the present context because of the source term in (\ref{form_g}).
\end{rk}

When $\omega$ is an arbitrary open subset of $\Omega$ and $\gamma \in (0,1)$,
we can still prove Lipschitz stability under an additional smallness assumption of the source term, which is probably due just to technical reasons.

\begin{thm} \label{Main-thm}
Let $\gamma \in (0,1)$, $\omega$ be an open subset of $\Omega$.
Then, for every $T_0 \in [0,T_1)$, there exists $\eta=\eta(T_0)>0$ such that
for every $R$ satisfying (\ref{Hyp:R}) and
\begin{equation} \label{Hyp_R_small_var}
\frac{1}{R_0} \left( \int_{T_0}^{T_1} \| \partial_t R(t)\|_{L^\infty(\Omega_1)}^2 dt\right)^{1/2} < \eta
\end{equation}
system (\ref{Grushin}) satisfies the Lipschitz stability estimate on $(T_0,T_1)\times\omega$.
\end{thm}






\subsubsection{Observability and null controllability}

In this article, we are also interested in the observability  problem for (\ref{Grushin}).

\begin{Def} [Observability] 
Let $T>0$. System (\ref{Grushin}) is {\em observable in $\omega$ in time} $T$ if
there exists $C>0$ such that, for every $u_{0} \in L^{2}(\Omega)$, the solution of
\begin{equation} \label{adjointL2}
\left\lbrace \begin{array}{ll}
\partial_t u - \Delta_x u - |x|^{2\gamma} b(x) \Delta_y u=0 & (t,x,y) \in (0,T)\times\Omega\, ,\\
u(t,x,y)=0 & (t,x,y) \in (0,T)\times \partial \Omega\, ,\\
u(0,x,y)=u_0(x,y) & (x,y) \in \Omega\, ,
\end{array}\right.
\end{equation}
satisfies
$$\int_{\Omega} |u(T,x,y)|^{2} dx dy 
\leqslant C \int_{0}^{T} \int_{\omega} |u(t,x,y)|^{2} dx dy dt\, .$$
\end{Def}

As a corollary of the analysis developed for the proof of Theorems \ref{Main-thm-strip} and \ref{Main-thm} (see Remarks~\ref{re:1} and \ref{re:2}),
we obtain a direct proof of observability for Grushin-type parabolic equations.
The following statement is a generalization to the multidimensional case of \cite[Theorem 2]{Grushin} (where $N_1=N_2=1$ is assumed).

\begin{thm} \label{thm:obs}
Let $\omega$ be an open subset of $\Omega$. 
\begin{enumerate}
\item If $\gamma \in (0,1)$, then system (\ref{Grushin}) is observable in $\omega$ in any time $T>0$.
\item If $\gamma=1$ and $\omega=\omega_1 \times \Omega_2$ where $\omega_1$ is an open subset of $\Omega_1$
then there exists $T^*>0$ such that for every $T>T^*$  system (\ref{Grushin}) is observable in $\omega$ in time $T$.
\end{enumerate}
\end{thm}

Note that we do not require $0 \in \omega$: the problem would be easily solved by cut-off functions arguments, under such an assumption (see \cite{Grushin}).
As a consequence, we deduce the following null controllability result.

\begin{Def}[Null controllability] 
Let $T>0$. System (\ref{Grushin}) is {\em null controllable in time} $T$ if,
for every $u_0 \in L^2(\Omega)$, there exists $g \in L^{2}((0,T)\times \Omega)$ such that the solution of
\begin{equation} \label{Grushin_CYpb}
\left\lbrace \begin{array}{ll}
\partial_t u - \Delta_x u - |x|^{2\gamma} b(x) \Delta_y u= g(t,x,y)1_{\omega}(x,y) &(t,x,y) \in (0,T)\times\Omega\,,\\
u(t,x,y)=0 & (t,x,y) \in (0,T)\times \partial \Omega\,,\\
u(0,x,y)=u_0(x,y) & (x,y) \in \Omega\,,
\end{array}\right.
\end{equation}
satisfies $u(T,\cdot,\cdot)=0$.
\end{Def}

Here $1_\omega$ is the characteristic function of the set $\omega$.

\begin{thm} \label{thm:cont}
Let $\omega$ be an open subset of $\Omega$.
\begin{enumerate}
\item If $\gamma \in (0,1)$, then system (\ref{Grushin}) is null controllable in any time $T>0$.
\item If $\gamma=1$ and $\omega=\omega_1 \times \Omega_2$ where $\omega_1$ is an open subset of $\Omega_1$, 
then there exists $T^*>0$ such that for every $T>T^*$ system (\ref{Grushin}) is null controllable in time $T$.
\end{enumerate}
\end{thm}

\subsection{Motivation and bibliographical comments}

We recall that the null controllability of (\ref{Grushin}), in the 2D case 
(i.e., $N_1=N_2=1$), is 
studied in detail in \cite{Grushin}. In particular, in 2D, null controllability:
\begin{itemize}
\item holds in any positive time $T>0$ with controls supported in an arbitrary open set $\omega$  when $\gamma \in (0,1)$,
\item holds only in large time $T>T^*>0$ when $\gamma=1$ and 
$\omega:= \omega_1 \times \Omega_2$ is a strip parallel to the $y$-axis,  not containing the line segment $x=0$, and 
\item does not hold when $\gamma>1$.
\end{itemize}
The goal of this article is:
\begin{itemize}
\item to generalize the previous positive controllability results to the multidimensional case, and
\item to prove a Lipschitz stability estimate for the inverse source problem,
by adapting a method by Imanuvilov and Yamamoto~\cite{IY1}, 
for the values of $\gamma$ for which  null controllability holds.
\end{itemize}

Our formulation of the inverse problem corresponds to a single measurement
(see also Bukhgeim and Klibanov~\cite{BK} which first proposed a methodology 
based on Carleman estimates).  Following \cite{BK},  
many works have been published on this subject.  For uniformly parabolic equations we can refer the reader, for example, to 
Imanuvilov and Yamamoto~\cite{IY1}, Isakov~\cite{Is},
Klibanov~\cite{Kl}, Yamamoto~\cite{Ya}, and the references therein (the present list of references is by no means complete).
As for inverse problems for degenerate parabolic equations, see Cannarsa, Tort and Yamamoto~\cite{CTY1,CTY2}.

\subsection{Structure of the article}

This article is organized as follows.
\\

Section \ref{sec:WP} is devoted to preliminary results concerning the well posedness of (\ref{Grushin}),
the Fourier decomposition of its solutions, and the dissipation speed of the Fourier modes.

In Section \ref{sec:Carl}, we state a Carleman estimate for a heat equation
with nonsmooth coefficients, solved by the Fourier modes of the solution of (\ref{Grushin}).

Section \ref{sec:InvPb} is devoted to the proof of the Lispchitz stability estimates,
for the inverse source problem, i.e.  Theorems \ref{Main-thm-strip} and \ref{Main-thm}.



In the appendix we prove the Carleman estimate stated in Section \ref{sec:Carl}.

\subsection{Notation}

The Euclidian norm in $\mathbb{R}^N$ is denoted by $|.|$ for every $N \in \mathbb{N}^*$.
The notation $\|.\|$ refers to $L^2$-norms in the space variables $x$, $y$ or $(x,y)$,
depending on the context.  $1_\omega$ is the caracteristic function of the set $\omega$.

\section{Preliminaries}
\label{sec:WP}

\subsection{Well posedness}

\begin{Prop} \label{Prop:WP_ut}
Let $\gamma \in (0,1]$, $u_0 \in D(G_\gamma)$, $g \in H^1(0,T;L^2(\Omega))$, and
$u \in C^0([0,T];L^2(\Omega)) \cap L^2(0,T;H^1_\gamma(\Omega))$ solution of (\ref{Grushin}) (\ref{IC}). 
Then the function $v:=\partial_t u$ belongs to $L^2(0,T;H^1_\gamma(\Omega))$ and solves
\begin{equation} 
\left\lbrace \begin{array}{ll}
\partial_t v - \Delta_x v - b(x) |x|^{2\gamma} \Delta_y v = \partial_t g(t,x,y)         & (t,x,y) \in (0, \infty)\times\Omega\,,\\
v(t,x,y)=0                                                                         & (t,x,y) \in (0, \infty)\times\partial \Omega\,,\\
v(0,x,y)=-G_\gamma u_0(x,y) + g(0,x,y)                                             & (x,y) \in \Omega\,.
\end{array}\right.
\end{equation}
\end{Prop}

\subsection{Fourier decomposition}
\label{subsec:Fourier}

We introduce the operator $\mathcal{A}$ defined by
$$D(\mathcal{A}):=H^2 \cap H^1_0(\Omega_2), \quad \mathcal{A}\varphi := - \Delta_y \varphi,$$
the nondecreasing sequence $(\mu_n)_{n \in \mathbb{N}^*}$ of its eigenvalues,
and the associated eigenvectors $(\varphi_n)_{n \in \mathbb{N}^*}$
\begin{equation} \label{def:varphin}
\left\lbrace \begin{array}{ll}
- \Delta_y \varphi_n(y) = \mu_n \varphi_n(y) & y \in \Omega_2,\\
\varphi_n(y)=0                               & y \in \partial \Omega_2.
\end{array}\right.
\end{equation}

\begin{Prop}
Let $u_0 \in L^2(\Omega)$, $g \in L^2((0,T)\times \Omega)$ and  $u$ be the solution of (\ref{Grushin}) (\ref{IC}).
For every $n \in \mathbb{N}^*$, the function
$$u_n(t,x):=\int_{\Omega_2} u(t,x,y) \varphi_n(y) dy$$
belongs to $C^0([0,T];L^2(\Omega))$ and is the unique weak solution of
\begin{equation} \label{Grushin_n}
\left\lbrace \begin{array}{ll}
\partial_t u_n - \Delta_x u_n + \mu_n |x|^{2\gamma} b(x) u_n = g_n(t,x)         & (t,x) \in (0,T)\times \Omega_1,\\
u_n(t,x)=0                                                                 & t \in (0,T)\times\partial\Omega_1,\\
u_n(0,x)=u_{n,0}(x)                                                        & x \in \Omega_1,
\end{array}\right.
\end{equation}
where
$$g_n(t,x):=\int_{\Omega_2} g(t,x,y) \varphi_n(y) dy
\quad \text{ and } \quad
u_{0,n}(x)=\int_{\Omega_2} u_0(x,y) \varphi_n(y) dy.$$
\end{Prop}

The proof is done as in \cite{Grushin}.

\subsection{Dissipation speed}
\label{subsec:dissipation}

We introduce, for every $n \in \mathbb{N}^*, \gamma>0$, the operator $G_{n,\gamma}$ defined 
on $L^2(\Omega_1)$ by
\begin{equation} \label{def:An}
\begin{array}{ll}
D(G_{n,\gamma}):=H^2 \cap H^1_0(\Omega_1)\,, 
&
G_{n,\gamma} u := -\Delta_x u + \mu_n |x|^{2\gamma} b(x) u .
\end{array}
\end{equation}
The smallest eigenvalue of $G_{n,\gamma}$ is given by
$$
\displaystyle\lambda_{n,\gamma} = \min \left\{
\frac{ \int_{\Omega_1} \left[ |\nabla v(x)|^2 + \mu_n |x|^{2\gamma} b(x) v(x)^2 \right] dx }{\int_{\Omega_1} v(x)^2 dx } ;
v \in H^1_0(\Omega_1) ,\ v\neq 0\right\}\, .
$$
We are interested in the asymptotic behavior (as $n \rightarrow + \infty$) of $\lambda_{n,\gamma}$,
which quantifies the dissipation speed of the solutions of (\ref{Grushin}).
The following result turns out to be a key point of the proofs of this article.

\begin{Prop} \label{Prop:1st_eigenvalue}
For every $\gamma>0$, there exists constants $c_*, c^* >0$ such that
$$ c_* \mu_n^{\frac{1}{1+\gamma}} \leqslant \lambda_{n,\gamma} 
\leqslant c^* \mu_n^{\frac{1}{1+\gamma}}, \qquad \forall\, n \in \mathbb{N}^*\,.$$
\end{Prop}

\noindent \textbf{Proof of Proposition \ref{Prop:1st_eigenvalue}:}
First, we prove the lower bound.
Let $\tau_n:=\mu_n^{\frac{1}{2(1+\gamma)}}$.
With the change of variable $\phi(x)=\tau_n^{N_1/2} \varphi(\tau_n x)$, we get
$$\begin{array}{l}
\lambda_{n,\gamma}
= \inf \left\{ \int\limits_{\Omega_1} \left( |\nabla \phi(x)|^2 + \mu_n |x|^{2\gamma} b(x) \phi(x)^2 \right) dx ; 
\phi \in C^\infty_c(\Omega_1), \|\phi\|_{L^2(\Omega_1)}=1 \right\}
\\
= \tau_n^2
\inf\left\{ \int\limits_{\tau_n \Omega_1}
\left( |\nabla \varphi(y)|^2 +  |y|^{2\gamma} b(y/\tau_n) \varphi(y)^2 \right) dy ;
\varphi \in C^\infty_c(\tau_n \Omega_1), 
\|\varphi\|_{L^2(\tau_n \Omega_1)}=1 \right\}
\\
\geqslant c_* \tau_n^2
\end{array}$$
where
$$c_*:=\inf\left\{ 
\int_{\mathbb{R}^{N_1}} \left( |\nabla \varphi(y)|^2 + |y|^{2\gamma} b_* \varphi(y)^2 \right) dy ;
\varphi \in C^\infty_c(\mathbb{R}^{N_1}), \|\varphi\|_{L^2(\mathbb{R}^{N_1})}=1 \right\}$$
is positive (see \cite{resi}) and $b_*:=\min\{ b(x) ; x \in \overline{\Omega_1}\}$.

Now, we prove the upper bound in Proposition \ref{Prop:1st_eigenvalue}.
For every $k>1$ let us consider the function
$\varphi_k(x):= (1 - k|x|)^+$, that belongs to $H^1_0(\Omega)$ for $k$ large enough (so that $B_{\mathbb{R}^{N_1}}(0,1/k) \subset \Omega_1$). 
Easy computations show that
$$\begin{array}{c}
\displaystyle \int_{\Omega_1} \varphi_k(x)^2 dx = C_1(N)k^{-N}\, ,\ 
\displaystyle\int_{\Omega_1} |\nabla \varphi_k(x)|^2 dx = C_N(N)k^{2-N}\, ,
\\
\displaystyle \int_{\Omega_1} \mu_n |x|^{2\gamma} b(x) \varphi_k(x)^2 dx \leqslant C_3(N,\gamma) \mu_n b^* k^{-N-2\gamma}\,.
\end{array}$$
where $b^*:=\max\{b(x);x\in \overline{\Omega_1} \}$.
Thus, 
$$\lambda_{n,\gamma} \leqslant f_{n,\gamma}(k):= [C_2 k^2 + C_3 \mu_n b^* k^{-2\gamma}]/C_1, \quad \forall k>1.$$
Minimizing the right-hand side over $k$, we get $\lambda_{n,\gamma} \leqslant C(N,\gamma) \mu_n^{\frac{1}{\gamma+1}}$. \hfill $\Box$

\section{Carleman inequality for heat equations with nonsmooth potentials}
\label{sec:Carl}

For $\mu > 0$, let us introduce the operator
$$\mathcal{P}_{\mu,\gamma} u := \frac{\partial u}{\partial t} - \Delta_x u + \mu |x|^{2\gamma} b(x) u.$$
The goal of this section is the statement of the following Carleman inequality.

\begin{Prop} \label{Carleman_global}
Let $\gamma \in (0,1]$.
There exist a weight function $\beta \in C^1(\Omega_1;(0,\infty))$ and
positive constants $\mathcal{C}_1, \mathcal{C}_2$ such that
for every $\mu \in (0,\infty)$, $0\leqslant T_0<T_1 \leqslant T$, and $u \in C^0([0,T];L^2(\Omega_1)) \cap L^2(0,T;H^1_0(\Omega_1))$
the following inequality holds
\begin{equation} \label{Carl_est}
\begin{array}{ll}
& \mathcal{C}_1 \int_{T_0}^{T_1} \int_{\Omega_1} \left(
\frac{M}{(t-T_0)(T_1-t)} \big| \nabla_x u (t,x) \big|^2 +
\frac{M^3}{((t-T_0)(T_1-t))^3} \big| u(t,x) \big|^2 
\right)  e^{-M \alpha(t,x)} dx dt
\\
\leqslant &
\int_{T_0}^{T_1} \int_{\Omega_1} | \mathcal{P}_{\mu,\gamma} u(t,x) |^2 e^{-M\alpha(t,x)} dx dt 
\\ & +
\int_{T_0}^{T_1} \int_{\omega_1} \frac{M^3}{((t-T_0)(T_1-t))^3} | u(t,x)|^2 e^{-M\alpha(t,x)} dx dt
\end{array}
\end{equation}
where 
$$\alpha(t,x):=\frac{\beta(x)}{(t-T_0)(T_1-t)},$$
$$M:=\left\lbrace \begin{array}{l}
\mathcal{C}_2 \max\{T+T^2;\sqrt{\mu} T^2\} \text{ if } \gamma \in [1/2,1], \\
\mathcal{C}_2 \max\{T+T^2;\mu^{2/3} T^2\}  \text{ if } \gamma \in (0,1/2),
\end{array}\right.$$
and $T:=T_1-T_0$.
\end{Prop}

Note that we can have sharp dependency of $M = O(\mu^{1/2})$ and $T$
in the case of $1/2\le \gamma \le 1$.
In particular, if we treat the term $\mu\vert x\vert^{2\gamma}b(x)u$ as
lower-order term to apply the Carleman estimate for the operator
$\frac{\partial}{\partial t} - \Delta_x$, then we can 
obtain less sharp dependency $M = O(\mu^{2/3})$ and we need sharper 
estimate for $1/2\le \gamma \le 1$.

The proof of this Carleman inequality is given in \cite{Grushin} 
in the case $N_1=1$. In the 1D case, the sharp dependency  $M = O(\mu^{1/2})$
is proved for any $\gamma \in (0,1)$ and the case $\gamma \in (0,1/2)$ requires 
a weight adapted to the degeneracy.

A proof in the multi-dimensional case is presented in Appendix. 
It relies on the usuall weight of heat equations.

\section{Inverse source problem}
\label{sec:InvPb}

\subsection{Uniform observability of frequencies}

\begin{Prop} \label{Prop:UO}
Let $\gamma \in (0,1)$ and $\omega_1$ be an open subset of $\Omega_1$.
There exists $C>0$ and functions $\epsilon_n: (0,+\infty) \rightarrow 
(0,+\infty)$, $n \in \mathbb{N}^*$ with
\begin{itemize}
\item $\epsilon_n(T) \rightarrow 0$ when $n \rightarrow 0$, for every $T>0$,
\item $\epsilon_n(T) \leqslant \epsilon^* < + \infty$, for every $n \in \mathbb{N}^*$ and $T>0$,
\end{itemize}
such that, for every $n \in \mathbb{N}^*$, $g_n \in L^2((0,T)\times\Omega_1)$, $u_{0,n} \in L^2(\Omega_1)$ the solution of (\ref{Grushin_n}) satisfies
\begin{equation} \label{Ineq:UO}
\begin{array}{ll}
\int_{\Omega_1} |u_n(T,x)|^2 dx
\leqslant 
  e^{C\left(1+T^{-p}\right)}
\Big( &
\int_0^T \int_{\omega_1} |u_n(t,x)|^2 dx dt  
\\ & 
+ \epsilon_n(T) \int_0^T \int_{\Omega_1} |g_n(t,x)|^2 dx dt
\Big)
\end{array}
\end{equation}
where
\begin{equation} \label{def:p}
p=p(\gamma):=\left\lbrace \begin{array}{l}
\frac{1+\gamma}{1-\gamma}, \text{ if } \gamma \in [1/2,1],\\
\frac{2(1+\gamma)}{1-2\gamma}, \text{ if } \gamma \in (0,1/2).
\end{array}\right.
\end{equation}
\end{Prop}

\noindent \textbf{Proof of Proposition \ref{Prop:UO}:} The proof is in 3 steps.
\\

\noindent \emph{Step 1: We prove}
\begin{equation} \label{step1}
\|u_n(T)\|^2_{L^2(\Omega_1)} \leqslant 
\frac{6}{T}  e^{-2 \lambda_{n,\gamma}T/3} \int_{T/3}^{2T/3} \|u_n(t)\|^2 dt + 
\frac{1}{\lambda_{n,\gamma}} \|g_n\|_{L^2((0,T)\times\Omega_1)}^2.
\end{equation}
From Duhamel's formula, i.e.,
$$u_n(T)=e^{-G_{n,\gamma}(T-t)}u_n(t)+\int_t^T e^{-G_{n,\gamma}(T-\tau)} g_n(\tau) d\tau, \quad \forall t \in (0,T)$$
and the Cauchy-Schwarz inequality, we get 
$$\begin{array}{ll}
\|u_n(T)\|_{L^2(\Omega_1)} 
& \leqslant e^{-\lambda_{n,\gamma}(T-t)} \|u_n(t)\|_{L^2(\Omega_1)} 
+ \int_t^T e^{-\lambda_{n,\gamma}(T-\tau)} \|g_n(\tau)\|_{L^2(\Omega_1)} d\tau
\\ & \leqslant
e^{-\lambda_{n,\gamma}(T-t)} \|u_n(t)\|_{L^2(\Omega_1)}  
+ \frac{1}{\sqrt{2\lambda_{n,\gamma}}} \|g_n\|_{L^2((0,T)\times\Omega_1)}.
\end{array}$$
Thus
$$\|u_n(T)\|^2_{L^2(\Omega_1)} \leqslant
2 e^{-2 \lambda_{n,\gamma}(T-t)} \|u_n(t)\|^2_{L^2(\Omega_1)} +
\frac{1}{\lambda_{n,\gamma}} \|g_n\|_{L^2((0,T)\times\Omega_1)}^2.$$
Integrating this relation over $t \in (T/3,2T/3)$ gives (\ref{step1}).
\\

\noindent \emph{Step 2: We prove the existence of a constant $\mathcal{C}_3>0$ such that
for every $T>0$, $n \in \mathbb{N}^*$, $g_n \in L^2((0,T)\times\Omega_1)$
$u_{0,n} \in L^2(\Omega_1)$, the solution of (\ref{Grushin_n}) satisfies}
\begin{equation} \label{step2}
\int_{T/3}^{2T/3} \int_{\Omega_1} |u_n|^2 dxdt 
\leqslant  \mathcal{C}_3 T e^{\frac{9M\beta^*}{2T^2}} \left( 
\int_0^T \int_{\Omega_1} |g_n|^2 + \int_0^T \int_{\omega_1} |u_n|^2  \right)
\end{equation}
where $\beta$, $\mathcal{C}_2$ and $M$ are as in Proposition \ref{Carleman_global} (with $\mu$ replaced by $\mu_n$), 
and $\beta^*:=\max\{\beta(x);x \in \Omega_1 \}$. 
From Proposition \ref{Carleman_global}, we get
\begin{equation} \label{step2_1}
\begin{array}{ll}
& \displaystyle \mathcal{C}_1 \left( \frac{4M}{T^2} \right)^3 e^{-\frac{9M\beta^*}{2T^2}} \int_{T/3}^{2T/3} \int_{\Omega_1} |u_n|^2 dxdt 
\\ \leqslant &
\displaystyle \mathcal{C}_1 \int_{T/3}^{2T/3} \int_{\Omega_1} \frac{M^3}{(t(T-t)])^3} |u_n|^2  e^{-\frac{M \beta(x)}{t(T-t)}} dxdt
\\ \leqslant &
\displaystyle \mathcal{C}_1 \int_{\Omega_1}\int^T_0 
\frac{M^3}{(t(T-t))^3} |u_n|^2 e^{-\frac{M\beta(x)}{t(T-t)}} dxdt \\ 
\leqslant &
\displaystyle \int_0^T \int_{\Omega_1} |g_n|^2 e^{-\frac{M \beta(x)}{t(T-t)}} 
dxdt + \int_0^T \int_{\omega_1} \frac{M^3}{(t(T-t)])^3} |u_n|^2  
e^{-\frac{M \beta(x)}{t(T-t)}} dxdt \\ 
\leqslant &
\displaystyle \int_0^T \int_{\Omega_1} |g_n|^2 dxdt 
+ C \int_0^T \int_{\omega_1} |u_n|^2 dx dt
\end{array}
\end{equation}
where $C:=\sup\{ x^3 e^{-\beta_* x} ; x \geqslant 0\}$ and $\beta_*:=\min\{\beta(x);x\in\omega_1\}$.
We deduce from (\ref{step2_1}) that
$$
\int_{T/3}^{2T/3} \int_{\Omega_1} |u_n(t)|^2 dxdt 
\leqslant \frac{\max\{1,C\}}{4^3 \mathcal{C}_1} \frac{T^6}{M^3} e^{\frac{9M\beta^*}{2T^2}} \left( 
\int_0^T \int_{\Omega_1} |g_n|^2  + \int_0^T \int_{\omega_1} |u_n|^2  \right).$$
We remark that $M \geqslant \mathcal{C}_2 T$ and $\mathcal{C}_2 T^2$ thus 
$T^6/M^3 \leqslant T/\mathcal{C}_2^3$. Then, the previous inequality
gives (\ref{step2}) with $\mathcal{C}_3:=\max\{1,C\}/(4^3 \mathcal{C}_1 \mathcal{C}_2^3)$.
\\

\noindent \emph{Step 3: We put together (\ref{step1}) and (\ref{step2}).} 
Using (\ref{step1}) in the first inequality, (\ref{step2}) in the second one, and Proposition \ref{Prop:1st_eigenvalue} in the third one, we get
\begin{equation} \label{Obs_interm}
\begin{array}{lll}
\displaystyle \int_{\Omega_1} |u_n(T)|^2
& \leqslant &
\displaystyle \frac{6}{T}  e^{-2 \lambda_{n,\gamma}T/3} \int_{T/3}^{2T/3} 
\|u_n(t)\|^2_{L^2(\Omega_1)} dt 
+ \frac{1}{\lambda_{n,\gamma}} \|g_n\|_{L^2((0,T)\times\Omega_1)}^2
\\ 
& \leqslant &
\displaystyle 6 \mathcal{C}_3  e^{\frac{9M\beta^*}{2T^2}-2 \lambda_{n,\gamma}T/3} \int_0^T \int_{\omega_1} |u_n|^2 dxdt
\\ & & 
\displaystyle + \left( \frac{1}{\lambda_{n,\gamma}} + 6 \mathcal{C}_3  e^{\frac{9M\beta^*}{2T^2}-2 \lambda_{n,\gamma}T/3} \right) 
\int_0^T \int_{\Omega_1} |g_n|^2 dxdt\\ 
& \leqslant &
\displaystyle 6 \mathcal{C}_3  e^{\frac{9M\beta^*}{2T^2}-c_1 T \mu_n^{\frac{1}{1+\gamma}}}  \int_0^T \int_{\omega_1} |u_n|^2 dxdt
\\ & &
\displaystyle +\left( \frac{1}{c_* \mu_n^{\frac{1}{1+\gamma}}} + 6 \mathcal{C}_3  e^{\frac{9M\beta^*}{2T^2}-c_1 T \mu_n^{\frac{1}{1+\gamma}}} \right) \int_0^T \int_{\Omega_1} |g_n|^2,
\end{array}
\end{equation}
where $c_1:=2c_*/3$.
\\

\noindent \emph{Step 4: End of the proof when $\gamma \in [1/2,1]$.}

\noindent \emph{\underline{First case: $\mu_n \geqslant 1+\frac{1}{T}$.}} Then we have that $M=\mathcal{C}_2 \sqrt{\mu_n} T^2$.
For any constant $c_2>0$, the maximum value of the function $z \mapsto c_2 z - c_1 z^{\frac{2}{1+\gamma}} T/2$ 
on $(0,+\infty)$ is of the form $c_3 T^{-\frac{1+\gamma}{1-\gamma}}$ for some constant $c_3>0$ (independent of $T$). Thus,
$$\begin{array}{ll}
\|u_n(T)\|^2_{L^2(\Omega_1)} \leqslant&   
6 \mathcal{C}_3 e^{c_3 T^{-\frac{1+\gamma}{1-\gamma}}}  
\int_0^T \int_{\omega_1} |u_n|^2 dxdt
\\ &
+\left( \frac{1}{c_* \mu_n^{\frac{1}{1+\gamma}}} + 6 \mathcal{C}_3 e^{c_3 T^{-\frac{1+\gamma}{1-\gamma}}-\frac{c_1 T}{2} \mu_n^{\frac{1}{1+\gamma}}} \right)
\int_0^T \int_{\Omega_1} |g_n|^2 dxdt.
\end{array}$$
This proves (\ref{Ineq:UO}) with any constant $C$ large enough so that
$$
6 \mathcal{C}_3 \leqslant e^C, \qquad \text{ and } \qquad
c_3 \leqslant C$$
and
$$\epsilon_n(T):=\frac{1}{c_* \mu_n^{\frac{1}{1+\gamma}}} + 6 \mathcal{C}_3 e^{-\frac{c_1 T}{2} \mu_n^{\frac{1}{1+\gamma}}}.$$

\noindent \emph{\underline{Second case: $\mu_n<1+\frac{1}{T}$.}} Then, $M=\mathcal{C}_2(T+T^2)$ thus
$$\begin{array}{ll}
\|u_n(T)\|^2_{L^2(\Omega_1)} \leqslant
&   6 \mathcal{C}_3  e^{\frac{9\beta^*}{2} \mathcal{C}_2 \left( 1 + \frac{1}{T} \right)}  \int_0^T \int_{\omega_1} |u_n|^2 dxdt\\ 
& +\left( \frac{1}{c_* \mu_0^{\frac{1}{1+\gamma}}} + 6 \mathcal{C}_3  e^{\frac{9\beta^*}{2} \mathcal{C}_2 \left( 1 + \frac{1}{T} \right)}
 \right) \int_0^T \int_{\Omega_1} |g_n|^2 dxdt.
\end{array}$$
Note that $\frac{1+\gamma}{1-\gamma}>1$ thus
$\frac{1}{T} \leqslant 1$ when $T \geqslant 1$ and
$\frac{1}{T} \leqslant T^{-\frac{1+\gamma}{1-\gamma}}$ when $T <1$;
in any case $\frac{1}{T} \leqslant 1 + T^{-\frac{1+\gamma}{1-\gamma}}$. Thus, we have
$$\int\limits_{\Omega_1} |u_n(T)|^2_ \leqslant
\left( \frac{1}{c_* \mu_0^{\frac{1}{1+\gamma}}} + 6 \mathcal{C}_3  e^{\frac{9\beta^*}{2} \mathcal{C}_2 \left( 2 + T^{-\frac{1+\gamma}{1-\gamma}} \right)} \right)
\left( \int\limits_0^T \int\limits_{\omega_1} |u_n|^2  + \int\limits_0^T \int\limits_{\Omega_1} |g_n|^2 \right).$$
This proves (\ref{Ineq:UO}) with $\epsilon_n=1$ and any constant $C$ large enough so that
$$\frac{1}{c_* \mu_0^{\frac{1}{1+\gamma}}} \leqslant \frac{1}{2} e^C, \qquad
6 \mathcal{C}_3  e^{9\beta^* \mathcal{C}_2}  \leqslant \frac{1}{2} e^C, \qquad \text{ and } \qquad
\frac{9\beta^* \mathcal{C}_2}{2}  \leqslant C. $$

\noindent \emph{Step 5: End of the proof when $\gamma \in (0,1/2)$.}
One can proceed as in Step 4 observing that the maximum value of the function
$z \mapsto c_2 z - c_1 z^{\frac{3}{2(1+\gamma)}}T/2$ on $(0,+\infty)$ is of the form
$c_3 T^{-\frac{2(1+\gamma)}{1-2\gamma}}$ for the first case, and
$\frac{2(1+\gamma)}{1-2\gamma}>1$ for the second one. $\hfill \Box$
\\

\begin{Prop} \label{Prop:UO_gamma=1}
Assume $\gamma=1$ and let $\omega_1$ be an open subset of $\Omega_1$.
Then there exists $T_*>0$ such that,
for every $T>T_*$, there exists $C>0$ and 
$(\epsilon_n(T))_{n \in \mathbb{N}^*} \in (\mathbb{R}^*_+)^{\mathbb{N}^*}$ with $\epsilon_n(T) \rightarrow 0$
such that, for every $n \in \mathbb{N}^*$, $g_n \in L^2((0,T)\times\Omega_1)$, $u_{0,n} \in L^2(\Omega_1)$, 
the solution of (\ref{Grushin_n}) satisfies
$$\int_{\Omega_1} |u_n(T,x)|^2 dx
\leqslant  C \int_0^T \int_{\omega_1} |u_n(t,x)|^2 dx dt  
+ \epsilon_n(T) \int_0^T \int_{\Omega_1} |g_n(t,x)|^2 dx dt.$$
\end{Prop}
\begin{rk}\label{re:1}
The above proposition, together with the Bessel-Parseval equality, proves statement 2 of Theorem~\ref{thm:obs}.
\end{rk}
\noindent \textbf{Proof of Proposition \ref{Prop:UO_gamma=1}:}
One can follow the lines of the previous proof until (\ref{Obs_interm}).
Then, when $\mu_n \geqslant 1+\frac{1}{T}$, we have $M=\mathcal{C}_2 \sqrt{\mu_n} T^2$. Thus
$$\begin{array}{ll}
\|u_n(T)\|^2_{L^2(\Omega_1)} \leqslant
&   
6 \mathcal{C}_3 e^{[c_2-c_1 T] \sqrt{\mu_n}}  \int_0^T 
\int_{\omega_1} |u_n|^2 dxdt \\ 
& +\left( \frac{1}{\lambda_{n,\gamma}} + \mathcal{C}_3  e^{[c_2-c_1 T] \sqrt{\mu_n}} \right) \int_0^T \int_{\Omega_1} |g_n|^2 dxdt
\end{array}$$
for some constant $c_2>0$. This gives the conclusion with $T_*=c_2/c_1$. \hfill$\Box$

\subsection{Lipschitz stability estimate when $\omega$ is a strip}

The goal of this section is the proof of Theorem \ref{Main-thm-strip}.
We focus on the uniform Lipschitz stability for systems (\ref{Grushin_n}).
We assume the source term $g_n$ in (\ref{Grushin_n}) takes the form $g_n(t,x)=f_n(x) R(t,x)$,
where $f_n \in L^2(\Omega_1)$ and $R \in C^0([0,T]\times\overline{\Omega_1})$.

\begin{Def}[Uniform Lipschitz stability]
Let $\omega_1$ be an open subset of $\Omega_1$, $T>0$ and $0\leqslant T_0 < T_1 \leqslant T$.
We say the system (\ref{Grushin_n}) satisfies a \emph{uniform Lipschitz stability estimate on $(T_0,T_1)\times\omega_1$} if,
for every $n \in \mathbb{N}^*$, $f_n \in L^2(\Omega_1)$, $u_{0,n} \in L^2(\Omega_1)$,
the solution of (\ref{Grushin_n}) satisfies
\begin{equation} \label{Lipschitz-ineq-n}
\int_{\Omega_1} |f_n(x)|^2 dx \leqslant
C \left(
\int_{T_0}^{T_1} \int_{\omega_1} |\partial_t u_n(t,x)|^2 dx dt +
\int_{\Omega_1} |G_{n,\gamma}u_n(T_1,x)|^2 dx 
\right).
\end{equation}
\end{Def}

Theorem \ref{Main-thm-strip} is a consequence of the following proposition and the Bessel-Parseval equality.

\begin{Prop} \label{Prop:Stab_n}
Assume (\ref{Hyp:R}) and let $\omega_1$ be an open subset of $\Omega_1$.  
\begin{enumerate}
\item If $\gamma \in (0,1)$ then, for every $T_0 \in [0,T_1)$, 
system (\ref{Grushin_n}) satisfies a uniform Lipschitz stability estimate on $(T_0,T_1)\times\omega_1$.
\item If $\gamma=1$, then there exists $T^*>0$ such that, for every $T_0 \in [0,T_1-T^*)$,
system (\ref{Grushin_n}) satisfies a uniform Lipschitz stability estimate on $(T_0,T_1)\times\omega_1$.
\end{enumerate}
\end{Prop}

\begin{rk}
The inequality (\ref{Lipschitz-ineq-n}) with a constant $C$ that may depend on $n$ is already known (see \cite{IY1}). 
The goal of this section is to prove that (\ref{Lipschitz-ineq-n}) holds with a constant $C$ which is independent of $n$.
\end{rk}

\noindent \textbf{Proof of Proposition \ref{Prop:Stab_n}:} 
In this proof, $T^*$ is as in Proposition \ref{Prop:UO_gamma=1} if $\gamma=1$ and $T^*:=0$ if $\gamma \in (0,1)$.
Assume $(T_1-T_0)>T_*$. It results from  (\ref{Hyp:R}) that
$$R_0 |f_n(x)|  \leqslant |R(T_1,x) f_n(x)|  \leqslant |\partial_t u_n(T_1,x)| + |G_{n,\gamma} u(T_1,x)|.$$
and
\begin{equation} \label{interm1}
\int_{\Omega_1} |f_n(x)|^2 dx 
\leqslant \frac{2}{R_0^2} \left(
\int_{\Omega_1} |\partial_t u_n(T_1,x)|^2 dx +
\int_{\Omega_1} |G_{n,\gamma}u_n(T_1,x)|^2 dx 
\right).
\end{equation}
By Propositions \ref{Prop:WP} and \ref{Prop:UO} or \ref{Prop:UO_gamma=1} (applied to $\partial_t u$), we get, for every $T_0 \in (0,T_1-T^*)$,
\begin{equation} \label{interm2}
\int_{\Omega_1} |\partial_t u_n(T_1)|^2 dx
\leqslant C \int_{T_0}^{T_1} \int_{\omega_1} |\partial_t u_n|^2 dxdt
+ \epsilon_n [T_1-T_0] \|\partial_t R \|_{\infty}^2 \int_{\Omega_1} |f_n|^2 dx
\end{equation}
with $\epsilon_n=\epsilon_n(T_1-T_0)$. There exists $n^*>0$ such that, for every $n \geqslant n^*$,
$(2/R_0^2) \epsilon_n [T_1-T_0] \|\partial_t R \|_{\infty}^2 <1/2$.
Using (\ref{interm1}) and (\ref{interm2}) we get, for $n \geqslant n^*$
$$\int_{\Omega_1} |f_n|^2 dx
\leqslant 
\frac{4 C}{R_0^2}  \int_{T_0}^{T_1} \int_{\omega_1} |\partial_t u_n|^2 dxdt 
+ \frac{4}{R_0^2}  \int_{\Omega_1} |G_{n,\gamma}u_n(T_1)|^2 dx.
$$
This ends the proof of Proposition \ref{Prop:Stab_n}. $\hfill \Box$

\subsection{Lipschitz stability estimate when $\omega$ is arbitrary}

The goal of this section is the proof of Theorem \ref{Main-thm}.
In all the section, $T>0$ and $\gamma \in (0,1)$ are fixed. 
For simplicity, we take $T_0=0$ and $T_1=T$. $\Vert \cdot\Vert$ denotes the
norm in $L^2(\Omega)$. 
\\

For $n \in \mathbb{N}^*$, $\varphi_n(y)$ is defined by (\ref{def:varphin}) and
$H_n:= L^2(\Omega_1) \otimes \varphi_n$ is a closed subspace of $L^2(\Omega)$.
For $j \in \mathbb{N}^*$, we define 
$$E_j:=\oplus_{\mu_n \leqslant 2^{2j}} H_n$$
and denote by $\Pi_j$ the orthogonal projection onto $E_j$.  Moreover, $Id$ stands for the identity operator on $L^2(\Omega)$.
\\

\begin{Prop} \label{Prop:obs_j}
Let $\omega$ be an open subset of $\Omega$.
Then there exists $C>0$ such that, for every $T>0$, $j \in \mathbb{N}^*$, 
$u_0 \in E_j$ and $g \in L^2(0,T;E_j)$,
the solution of (\ref{Grushin}) satisfies
$$\begin{array}{ll}
\int_{\Omega} |u(T,x,y)|^2 dx dy
\leqslant 
&
 e^{C\left(2^j +T^{-p} \right)} \int_0^T \int_\omega |u(t,x,y)|^2 dx dy dt 
\\ & 
+ e^{C\left(1 +T^{-p} \right)} \int_0^T \int_{\Omega} |g(t,x,y)|^2 dx dy dt,
\end{array}$$
where $p=p(\gamma)$ is defined by (\ref{def:p}).
\end{Prop}

For the proof of  Proposition \ref{Prop:obs_j} we shall need the following inequality obtained in \cite{Lebeau-Robbiano} 
(see also \cite{Lebeau-LeRousseau}).

\begin{Prop} \label{Prop:LR}
Let $\omega_2$ be an open subset of $\Omega_2$.
There exists $C >0$ such that, for every $(b_k)_{k \in \mathbb{N}^*} \subset \mathbb{R}$ and $\mu>0$,
$$\sum\limits_{\mu_k \leqslant \mu} |b_k|^2 \leqslant C e^{C \sqrt{\mu}} 
\int_{\omega_2} \left| \sum\limits_{\mu_k \leqslant \mu} b_k \varphi_k(y) \right|^2 dy.$$
\end{Prop}

\noindent \textbf{Proof of Proposition \ref{Prop:obs_j}:} 
Let $\omega_j$ be an open subset of $\Omega_j$ for $j=1,2$ such that
$\omega_1 \times \omega_2 \subset \omega$.
Using Proposition \ref{Prop:UO} and the orthonormality of the functions $(\varphi_n)$ in $L^2(\Omega_2)$, we get
$$\begin{array}{ll}
&\displaystyle \int_{\Omega} |u(T,x,y)|^2 dx dy =  \sum_{\mu_n \leqslant 2^{2j}} \int_{\Omega_1} |u_n(T,x)|^2 dx
\\ 
\leqslant & \displaystyle e^{C\left(1+T^{-p} \right)} \sum_{\mu_n \leqslant 2^{2j}} 
\left( \int_0^T \int_{\omega_1} |u_n(t,x)|^2 dx dt + \epsilon^* \int_0^T \int_{\Omega_1} |g_n(t,x)|^2 dx dt  \right)
\\  \leqslant & \displaystyle
e^{C\left(2^j+T^{-p} \right)}
\int_0^T \int_{\omega_1} \int_{\omega_2} \left| \sum_{\mu_n \leqslant 2^{2j}} u_n(t,x) \varphi_n(y) \right|^2 dy dx dt
\\ & \displaystyle + e^{C\left(1+T^{-p} \right)} \int_0^T \int_{\Omega} |g(t,x,y)|^2 dx dy dt. 
\end{array}$$
where the constant $C$ may change from line to line. $\hfill\Box$
\\
\\

Let $\rho \in \mathbb{R}$ be such that
\begin{equation} \label{def:rho}
0 < \rho < \min\left\{ \frac{1-\gamma}{1+\gamma} , \frac{1}{p(\gamma)} \right\}
\end{equation} 
where $p(\gamma)$ is defined by (\ref{def:p}) and $K=K(\rho)>0$ is such that $2K\sum_{j=1}^\infty 2^{-j\rho}=T$.
For $j \in \mathbb{N}^*$, we introduce 
\begin{equation} \label{def:tauj}
\tau_j:=K 2^{-j\rho}
\end{equation}
and $\alpha_j:=\sum_{k=1}^j 2\tau_k$.
By convention $\alpha_0:=0$. Let 
\begin{equation} \label{def:InJn}
I_n:=(T-\alpha_{n-1}-\tau_n,T-\alpha_{n-1})
\quad \text{ and }
J_{n}:=(T-\alpha_n, T-\alpha_{n-1}), \quad 
\forall n \in \mathbb{N}^*.
\end{equation}

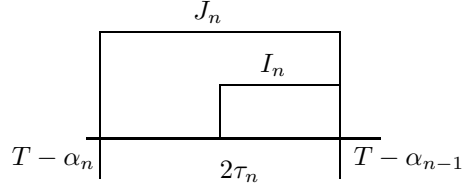
\bigskip
\begin{figure}[h]
\begin{picture}(450,40)
\put(215,-10){\line(0,1){55}} 
\put(125,45){\line(1,0){90}}
\put(170,25){\line(1,0){45}}
\put(170,25){\line(0,-1){20}}
\put(92,-5){$T-\alpha_n$}
\put(220,-5){$T-\alpha_{n-1}$} 
\put(160,50){$J_n$}
\put(185,30){$I_n$}
\put(170,-10){$2\tau_n$}
\put(125,-10){\line(0,1){55}}
\thicklines 
\put(120,5){\line(1,0){110}}
\end{picture}
\end{figure}

We will also use the notation 
\begin{equation} \label{lambda(2n)}
\lambda(2^n):=c_* 2^{\frac{2n}{1+\gamma}}
\end{equation} 
where $c_*$ is as in Proposition \ref{Prop:1st_eigenvalue} and we will write $G$ instead of $G_\gamma$.

\begin{Prop} \label{Prop:obs_n}
There exist $C_1, C_2, C_3>0$ such that for every
$n \in \mathbb{N}^*$, $g \in L^2((0,T)\times \Omega)$, $u_0 \in L^2(\Omega)$, 
the solution of (\ref{Grushin})(\ref{IC}) satisfies
\begin{equation}
\begin{array}{ll}
e^{-C_1 2^n} \| \Pi_n u (T-\alpha_{n-1}) \|^2 \leqslant
&
\int\limits_{I_n \times \omega} |u|^2 
+ C_2 2^{-n\left( \rho + \frac{2}{1+\gamma} \right)}  \int\limits_{J_n \times \Omega} |g|^2  \\
& + C_3 e^{-\lambda(2^n)\tau_n} \| (Id-\Pi_n)u(T-\alpha_n) \|^2.
\end{array}
\end{equation}
\end{Prop}
 
\noindent \textbf{Proof of Proposition \ref{Prop:obs_n}:} 
Let $n \in \mathbb{N}^*$, $g \in L^2((0,T)\times \Omega)$ and $u_0 \in L^2(\Omega)$.
By Proposition \ref{Prop:obs_j}, the solution of (\ref{Grushin})(\ref{IC}) satisfies
\begin{equation} \label{obs_n_interm0}
\begin{array}{ll}
\| \Pi_n u(T-\alpha_{n-1}) \|^2 \leqslant
& 
e^{C\left(2^n + \tau_n^{-p} \right)} \int\limits_{I_n \times \omega} | \Pi_n u |^2   \\
& 
+ e^{C\left(1 + \tau_n^{-p} \right)} \int\limits_{I_n \times \Omega} | \Pi_n g |^2  .
\end{array}
\end{equation}
Moreover, we have
\begin{equation} \label{obs_n_interm1}
\int\limits_{I_n \times \omega} |\Pi_n u|^2 
\leqslant \int\limits_{I_n \times \omega} 2 |u|^2  +  \int\limits_{I_n \times \Omega} 2 |(Id-\Pi_n)u|^2 .
\end{equation}
For every $t \in I_n$, the Duhamel formula, i.e.
$$\begin{array}{ll}
(Id-\Pi_n)u(t)= & e^{-G_\gamma(t-T+\alpha_n)}(Id-\Pi_n)u(T-\alpha_n)
\\ & +  \int_{T-\alpha_n}^t e^{-G_\gamma(t-\tau)} (Id-\Pi_n)g(\tau) d\tau,
\end{array}$$
Proposition \ref{Prop:1st_eigenvalue} and the Cauchy-Schwarz inequality give
$$\begin{array}{lll}
\| (Id-\Pi_n)u(t) \| 
& \leqslant &
\|(Id-\Pi_n)u(T-\alpha_n)\| e^{-\lambda(2^n) (t-T+\alpha_n)} 
\\ & & 
+ \int_{T-\alpha_n}^{t}  e^{-\lambda(2^n)(t-\tau)} \|(Id-\Pi_n)g(\tau)\| d\tau 
\\ &  \leqslant & 
\|(Id-\Pi_n)u(T-\alpha_n)\| e^{-\lambda(2^n) \tau_n} 
\\ & & + \frac{1}{\sqrt{2\lambda(2^n)}} \|(Id-\Pi_n)g\|_{L^2(J_n\times\Omega)}.
\end{array}$$
Thus
\begin{equation} \label{obs_n_interm2}
\begin{array}{ll}
\int\limits_{I_n \times \Omega} |(Id-\Pi_n)u|^2 
\leqslant
& 
2 \tau_n \|(Id-\Pi_n)u(T-\alpha_n)\|^2 e^{-2\lambda(2^n) \tau_n} 
\\ & 
+ \frac{\tau_n}{\lambda(2^n)} \|(Id-\Pi_n)g\|_{L^2(J_n\times\Omega)}^2.
\end{array}
\end{equation}
Using (\ref{obs_n_interm0}), (\ref{obs_n_interm1}) and (\ref{obs_n_interm2}), we get
$$\begin{array}{ll}
\| \Pi_n u(T-\alpha_{n-1}) \|^2 \leqslant
& 
2  e^{C\left(2^n + \tau_n^{-p}\right)} \int\limits_{I_n \times \omega} 
|u|^2
\\ &
+ \left( \frac{2\tau_n}{\lambda(2^n)} + e^{C(1-2^n)} \right)  e^{C\left(2^n + \tau_n^{-p}\right)}
\int_{J_n \times \Omega} |g|^2 
\\ &
+ 4 \tau_n e^{C\left(2^n + \tau_n^{-p} \right)-2\lambda(2^n) \tau_n} \|(Id-\Pi_n)u(T-\alpha_n)\|^2.
\end{array}$$
In view of (\ref{def:rho}), we have that
$$\tau_n^{-p}=K^{-p} 2^{n \rho p } \leqslant C 2^n, \quad \forall n \in \mathbb{N}^*,$$
for some constant $C>0$. Thus, there exists $C_1>0$ such that
$$2 e^{C\left(2^n + \tau_n^{-p} \right)} 
\leqslant e^{C_1 2^n}, \quad \forall n \in \mathbb{N}^*.$$
Using (\ref{def:tauj}), (\ref{lambda(2n)}) and (\ref{def:rho}), we obtain, for some constants $C_2, C_3>0$,
$$\frac{1}{2} \left( \frac{2\tau_n}{\lambda(2^n)} + e^{C(1-2^n)} \right) 
\leqslant
C_2 2^{-n\left( \rho + \frac{2}{1+\gamma} \right)}, \quad \forall n \in \mathbb{N}^*$$
and
$$ \tau_n e^{-2\lambda(2^n) \tau_n}
\leqslant
C_3 e^{-\lambda(2^n)\tau_n}, \quad \forall n \in \mathbb{N}^*.$$
Therefore, we have that
$$\begin{array}{ll}
\| \Pi_n u(T-\alpha_{n-1}) \|^2 \leqslant
& 
e^{C_1 2^n} \int\limits_{I_n \times \omega} |u|^2 
\\ &
+ C_2 2^{-n\left( \rho + \frac{2}{1+\gamma} \right)} e^{C_1 2^n}
\int_{J_n \times \Omega} |g|^2 
\\ &
+ C_3 e^{C_1 2^n -\lambda(2^n)\tau_n} \|(Id-\Pi_n)u(T-\alpha_n)\|^2,
\end{array}$$
for every $n \geqslant 1$.This gives the conclusion.\hfill$\Box$

\begin{Prop} \label{Prop:rec}
Let $T>0$.
Then, there exists $C=C(T)>0$ such that,
for every $u_0 \in L^2(\Omega)$, $g \in L^2((0,T)\times \Omega)$ of the form $g(t,x,y)=R(t,x)f(x,y)$ with $R \in L^\infty((0,T)\times\Omega_1)$ and $f \in L^2(\Omega)$,
the solution of (\ref{Grushin})-(\ref{IC}) satisfies
$$\int_\Omega |u(T,x,y)|^2 dx dy \leqslant
C\left(  \int_0^T \int_{\omega} |u|^2  + \int_0^T \|R(t)\|_{L^\infty(\Omega_1)}^2 dt \int\limits_{\Omega} |f|^2  \right).$$
\end{Prop}
\begin{rk}\label{re:2}
The above proposition can be used to prove statement 1 of Theorem~\ref{thm:obs}.
\end{rk}

\noindent \textbf{Proof of Proposition \ref{Prop:rec}:} Let $C_1$, $C_2$, $C_3$ be as in Proposition \ref{Prop:obs_n}.
\\

\noindent \emph{Step 1: We prove by induction on $n \in \mathbb{N}^*$ that, for every $n \in \mathbb{N}^*$,
$$(P_n): \quad \begin{array}{ll} 
\sum\limits_{k=1}^{n}  e^{-C_1 2^k} \|\Pi_k u(T-\alpha_{k-1})\|^2 \leqslant 
& \sum\limits_{k=1}^{n} \delta_k \int\limits_{I_k \times \omega} |u|^2
\\ & + A_n \left( \int_{T-\alpha_n}^{T} \| R(t)\|_{L^\infty(\Omega_1)}^2 dt \right) \|f\|^2 
\\ & + B_n \|u(T-\alpha_n)\|^2
\end{array}$$
where
\begin{equation} \label{ABe_n1}
\delta_1:=1, \quad
A_1:= C_2 2^{-\left( \rho + \frac{2}{1+\gamma} \right)}, \quad 
B_1:=C_3 e^{-\lambda(2^1)\tau_1}
\end{equation}
and 
\begin{equation} \label{rec:delta}
\delta_{n+1}:=\max\{2;1+B_n e^{C_1 2^n}\},
\end{equation}
\begin{equation} \label{rec:A}
A_{n+1}:= \max \left\{ A _n ;  \frac{B_n}{\lambda(2^{n+1})} + \delta_{n+1} C_2 2^{-(n+1)\left(\rho+\frac{2}{1+\gamma}\right)} \right\},
\end{equation}
\begin{equation} \label{rec:B}
B_{n+1}:= 2 B_{n} e^{-4\lambda(2^{n+1})\tau_{n+1}} + \delta_{n+1} C_3 e^{-\lambda(2^{n+1})\tau_{n+1}}.
\end{equation}
}

The inequality $(P_{1})$ is given by Proposition \ref{Prop:obs_n} with $n=1$. Indeed,
$$\begin{array}{ll}
\int\limits_{J_1 \times \Omega} |g|^2 
& = \int\limits_{J_1 \times \Omega} |R(t,x)f(x,y)|^2 dx dy dt
\\ & \leqslant 
 \int\limits_{J_1 \times \Omega} \|R(t)\|_{L^\infty(\Omega_1)}^2 |f(x,y)|^2 dx dy dt
\\ & \leqslant
\left(\int_{T-\alpha_1}^{T} \|R(t)\|_{L^\infty(\Omega_1)}^2 dt \right) \|f\|^2.
\end{array}$$

Let us now assume that $(P_n)$ holds for some $n \in \mathbb{N}^*$ and prove $(P_{n+1})$.
We have
$$B_n \|u(T-\alpha_n)\|^2 =  B_n \|\Pi_{n+1} u(T-\alpha_n)\|^2 +  B_n \|(Id-\Pi_{n+1}) u(T-\alpha_n)\|^2.$$
Moreover, using Duhamel's formula as in the previous proof, we get
$$\begin{array}{ll}
& \displaystyle     \|(Id-\Pi_{n+1}) u(T-\alpha_n)\|
\\ \leqslant &
   \displaystyle    \|(Id-\Pi_{n+1}) u(T-\alpha_{n+1})\| e^{-2 \lambda(2^{n+1}) \tau_{n+1}} 
\\ &  \displaystyle  + \int_{T-\alpha_{n+1}}^{T-\alpha_n} e^{-\lambda(2^{n+1}) (T-\alpha_n-s)} \| (Id-\Pi_{n+1}) R(s)f \| ds 
\\ \leqslant &
      \displaystyle  \|(Id-\Pi_{n+1}) u(T-\alpha_{n+1})\| e^{-2 \lambda(2^{n+1}) \tau_{n+1}} 
\\ &  \displaystyle  + \left( \int_{T-\alpha_{n+1}}^{T-\alpha_n} e^{-2 \lambda(2^{n+1}) (T-\alpha_n-s)} ds \right)^{1/2} 
       \left( \int_{T-\alpha_{n+1}}^{T-\alpha_n}  \| R(s) \|_{L^\infty(\Omega_1)}^2 ds \right)^{1/2} \| f \|
\\ \leqslant &
\displaystyle  \|(Id-\Pi_{n+1}) u(T-\alpha_{n+1})\| e^{-2 \lambda(2^{n+1}) \tau_{n+1}} 
\\ & \displaystyle  + \frac{1}{\sqrt{2 \lambda(2^{n+1})}} \left( \int_{J_{n+1}}  \| R(s) \|_{L^\infty(\Omega_1)}^2 ds \right)^{1/2} \| f \|.
\end{array}$$
Therefore,
$$\begin{array}{ll}
B_n \|u(T-\alpha_n)\|^2  \leqslant 
   &  B_n \|\Pi_{n+1} u(T-\alpha_n)\|^2 
\\ & + 2 B_n \|(Id-\Pi_{n+1}) u(T-\alpha_{n+1})\|^2 e^{-4 \lambda(2^{n+1}) \tau_{n+1}}
\\ & + \frac{B_n}{\lambda(2^{n+1})} \left( \int_{J_{n+1}}  \| R(s) \|_{L^\infty(\Omega_1)}^2 ds \right) \|f\|^2.
\end{array}$$
Thus, $(P_n)$ yields
\begin{equation} \label{rec:interm1}
\begin{array}{ll} 
          &  \sum\limits_{k=1}^{n}  e^{-C_1 2^k} \|\Pi_k u(T-\alpha_{k-1})\|^2 - B_{n} \|\Pi_{n+1} u(T-\alpha_n)\|^2 \\
\leqslant &  \sum\limits_{k=1}^{n} \delta_k \int\limits_{I_k \times \omega} |u|^2 \\
          &  + \left( A_n \int_{T-\alpha_n}^T \|R(t)\|_{L^\infty(\Omega_1)}^2 dt + \frac{B_n}{\lambda(2^{n+1})} \int_{J_{n+1}}  \| R(t) \|_{L^\infty(\Omega_1)}^2 dt \right) \|f\|^2 \\
          &  + 2 B_n e^{-4 \lambda(2^{n+1}) \tau_{n+1}} \|u(T-\alpha_{n+1})\|^2.
\end{array}
\end{equation}
Moreover, by Proposition \ref{Prop:obs_n}, we also have
\begin{equation} \label{rec:interm2}
\begin{array}{ll}
e^{-C_1 2^{n+1}} \| \Pi_{n+1} u (T-\alpha_{n}) \|^2 \leqslant
&
\int\limits_{I_{n+1} \times \omega} |u|^2  
\\ & + C_2  2^{- (n+1)\left( \rho + \frac{2}{1+\gamma} \right)} \left( \int_{J_{n+1}} \|R(t)\|_{L^\infty(\Omega_1)}^2 dt\right) \|f\|^2  
\\ & + C_3 e^{-\lambda(2^{n+1})\tau_{n+1}} \| u(T-\alpha_{n+1}) \|^2.
\end{array}
\end{equation}
Note that $\delta_{n+1}$ is chosen so that
$$\delta_{n+1} e^{-C_1 2^{n+1}} - B_n \geqslant e^{-C_1 2^{n+1}}.$$
Thus, summing (\ref{rec:interm1}) and $\delta_{n+1}*$(\ref{rec:interm2}), we get $(P_{n+1})$. This ends the first step.
\\

\noindent \emph{Step 2: We prove that $\delta_n = 2$ for $n$ large enough.}
Let $\tilde{B}_n:=B_ne^{C_1 2^n}$. For every $n \in \mathbb{N}^*$ we have either
$$\tilde{B}_{n+1}=2 \tilde{B}_{n} e^{-4\lambda(2^{n+1})\tau_{n+1}+C_1 2^n} + 2 C_3 e^{-\lambda(2^{n+1})\tau_{n+1}+C_1 2^{n+1}} $$
if $\delta_{n+1}=2$, or
$$\tilde{B}_{n+1}= 2 \tilde{B}_{n} e^{-4\lambda(2^{n+1})\tau_{n+1}+C_1 2^n} + (1+\tilde{B}_n) C_3 e^{-\lambda(2^{n+1})\tau_{n+1}+C_1 2^{n+1}}$$
if $\delta_{n+1}=1+B_n e^{C_1 2^n}$. 
Using Proposition \ref{Prop:1st_eigenvalue}, (\ref{def:tauj}) and
the inequality $\frac{2}{1+\gamma}-\rho>1$ (see (\ref{def:rho})), we get a constant $C>0$ such that
$\tilde{B}_{n+1} \leqslant C(\tilde{B}_n+C), \forall n \in \mathbb{N}^*$. We deduce the existence of another constant $C>0$ such that
$\tilde{B}_n \leqslant C^n$ for every $n \in \mathbb{N}^*$. Then, we have either
$$\tilde{B}_{n+1}\leqslant 2 C^n e^{-4\lambda(2^{n+1})\tau_{n+1}+C_1 2^n} + 2 C_3 e^{-\lambda(2^{n+1})\tau_{n+1}+C_1 2^{n+1}} $$
if $\delta_{n+1}=2$, or
$$\tilde{B}_{n+1}\leqslant 2 C^n e^{-4\lambda(2^{n+1})\tau_{n+1}+C_1 2^n} + (1+C^n) C_3 e^{-\lambda(2^{n+1})\tau_{n+1}+ C_1 2^{n+1}}$$
if $\delta_{n+1}=1+B_n e^{C_1 2^n}$. In any case, $\tilde{B}_n \rightarrow 0$ when $n \rightarrow \infty$ because
$\frac{2}{1+\gamma}-\rho>1$. Thus $\delta_n = 2$ for $n$ large enough.
\\

\noindent \emph{Step 3: We prove that $(A_n)$ is bounded.} 
By definition $(A_n)_{n \in \mathbb{N}^*}$ is a non decreasing sequence.
Moroever, 
$$ \frac{B_n}{\lambda(2^{n+1})} + \delta_{n+1} C_2 2^{-(n+1)\left(\rho+\frac{2}{1+\gamma}\right)} \underset{n \rightarrow + \infty}{\longrightarrow} 0.$$
Thus, for $n$ large enough, we have that
$$\frac{B_n}{\lambda(2^{n+1})} + \delta_{n+1} C_2 2^{-(n+1)\left(\rho+\frac{2}{1+\gamma}\right)} \leqslant A_1 \leqslant A_n.$$
This implies $A_{n+1}=A_n$ for all $n \in \mathbb{N}^*$.
\\

\noindent \emph{Step 4: We pass to the limit as $n \rightarrow \infty$ in $(P_n)$.}
The last term on the right-hand side of $(P_n)$ converges to zero because  $B_n \leqslant \delta^* e^{-C_1 2^n}$.
Thus, we get 
\begin{equation} \label{sum_infinie}
\begin{array}{ll}
\sum\limits_{n=1}^{\infty}  e^{-C_1 2^n}  \| \Pi_n u(T-\alpha_{n-1}) \|^2
\leqslant
&
C\left(  \int_0^T \int_{\omega} |u|^2 
+ \int_0^T \|R(t)\|_{L^\infty(\Omega_1)}^2 dt \int\limits_{\Omega} |f|^2 \right)
\end{array}
\end{equation}
with
$$C:=\max\Big\{\delta^* ; \max\{ A_n ; n \in \mathbb{N}^* \} \Big\}.$$

\noindent \emph{Step 5: Conclusion.}
Using Duhamel's formula and the convention $\Pi_0=0$, we get
$$\begin{array}{lll}
\| u(T) \|^2
& =         &\sum_{n=1}^{\infty} \|(\Pi_n-\Pi_{n-1})u(T)\|^2 \\
& \leqslant & \sum_{n=1}^\infty  2 \|(\Pi_n-\Pi_{n-1})u(T-\alpha_{n-1})\|^2 e^{-2\lambda(2^{n-1})\alpha_{n-1}} \\ 
&           & + \sum_{n=1}^\infty \frac{1}{\lambda(2^{n-1})} \left( \int_{T-\alpha_{n-1}}^{T} \|R(t)\|_{L^\infty(\Omega_1)}^2 dt \right)  \| (\Pi_n-\Pi_{n-1})f\|^2.
\end{array}$$
Then, by Proposition \ref{Prop:1st_eigenvalue}, (\ref{def:tauj}) and
the inequality $\frac{2}{1+\gamma}-\rho>1$, we obtain, for some constant $C>0$
\begin{equation} \label{stab2}
\|  u(T) \|^2
\leqslant C \sum\limits_{n=1}^{\infty} e^{-C_1 2^n} \| \Pi_n u(T-\alpha_{n-1}) \|^2
+ C \left( \int_{0}^{T} \|R(t)\|_{L^\infty(\Omega_1)}^2 dt \right)   \|f\|^2.
\end{equation}
Thus, (\ref{sum_infinie}) implies the conclusion. $\hfill \Box$
\\

\noindent \textbf{Proof of Theorem \ref{Main-thm}:} Using the equation in \eqref{Grushin},  and (\ref{Hyp:R}), we obtain
\begin{equation} \label{stab1}
\begin{array}{ll}
\Vert f\Vert^2
& \leqslant \frac{1}{R_0^2} \int_{\Omega} |R(T_1,x)f(x,y)|^2 dx dy \\
& \leqslant \frac{2}{R_0^2} \int_\Omega \Big( |\partial_t u(T_1,x,y)|^2 + |G_\gamma u(T_1,x,y)|^2 \Big) dx dy.
\end{array}
\end{equation}
Applying Propositions \ref{Prop:WP_ut} and \ref{Prop:rec}, we obtain, for some constant $C>0$, 
$$\| \partial_t u(T_1) \|^2 \leqslant C \left(\int_0^T \int_\omega |u|^2 dxdydt
+ \|\partial_t R\|_{L^2(0,T;L^\infty(\Omega_1))}^2 \|f\|^2 \right).$$
From (\ref{stab1}) and (\ref{stab2}), we get another constant $C>0$ such that
$$\|f\|^2 \leqslant \frac{C}{R_0^2} \left( \int_0^T \int_\omega |u|^2 
+ \|\partial_t R\|_{L^2(0,T;L^\infty(\Omega_1))}^2 \|f\|^2 + \|G_\gamma u(T_1)\|^2 \right).$$ 
We get the Lipschitz stability estimate if $\sqrt{C} \|\partial_t R\|_{L^2(0,T;L^\infty(\Omega_1))} /R_0 <1$. $\hfill \Box$

\begin{rk} \label{RK:Rt_pt}
The smallness assumption on $\partial_t R$ is used only to absorb the source term of the right-hand side
by the left-hand side, in the previous estimates. In other references (for example \cite{CTY1}),
the parameter $M$ in the Carleman estimate is chosen large enough for this absorption to be possible without
additionnal smallness assumptions on $R$. In our situation, to use the same trick, one would need dissipation
estimates in weighted $L^2$-spaces (weight given by the Carleman estimate), which may be quite difficult to prove.	

\end{rk}




\appendix

\section{Proof of the Carleman estimate}

The goal of this section is the proof of  Proposition \ref{Carleman_global}.
To simplify notations, we take $T_0=0$ and $T_1=T$.
Subsection \ref{subsec:weight} is devoted to the properties of an appropriate weight function;
the Carleman estimate in proved in Section \ref{subsec:Carlm_proof}.

\subsection{Properties of the weight function}
\label{subsec:weight}

\begin{Prop} \label{Prop:weight}
There exists $a, C_1, C_3 >0$ and $\beta \in 
C^4(\overline{\Omega_1}; (0,\infty))$ such that
\begin{equation} \label{beta_bord}
\frac{\partial \beta}{\partial \nu} \geqslant 0 \text{ on } \partial \Omega,
\end{equation}
\begin{equation} \label{beta_13}
\begin{array}{c}
(1-a) (\Delta \beta)(x) |Z|^2 - 2 D^2 \beta(x)(Z,Z) \geqslant C_1 |Z|^2\, ,  
\forall Z \in \mathbb{R}^{N_1}, x \in \Omega_1\setminus\widetilde{\omega}_1,
\\
(a-1) (\Delta \beta)(x) |\nabla \beta(x)|^2 - 2 D^2 \beta(x)(\nabla \beta(x),\nabla \beta(x)) \geqslant C_3\, ,
\forall x \in \Omega_1\setminus\widetilde{\omega}_1\,.
\end{array}
\end{equation}
\end{Prop}

Herafter, for any  $N \times N$ matrix $A$ and $Z \in \mathbb{R}^N$, we denote by $A(Z,Z)$ the scalar $Z^{T}AZ$.
\\

\noindent \textbf{Proof of Proposition \ref{Prop:weight}:} Let $a \in (1,3)$ and let 
$\psi \in C^4(\overline{\Omega_1})$ be such that
\begin{equation} \label{def:psi0}
\psi > 0 \text{ on } \overline{\Omega_1}, \quad
\psi = 0 \text{ on } \partial \Omega_1 \quad \text{ and } \quad
|\nabla\psi(x)|>0, \forall x \in \overline{\Omega_1}\setminus\widetilde{\omega}_1 
\end{equation}
(see \cite{Fursikov-Imanuvilov-186} or \cite[Lemma 2.68 on page 80]{JMC-book} for the existence of such a function).
Note that the $C^4$-regularity of the boundary of $\Omega_1$ ensures the $C^4$-regularity of the 
distance to the boundary of $\Omega_1$, which in turn allows to construct a $C^4$-function $\psi$ with the same construction
as in \cite[Lemma 2.68 on page 80]{JMC-book}.

There exist numbers $m_*,m^*>0$ such that
$$|\nabla\psi(x)|>m_*, \quad |\nabla \psi(x)|, |\Delta\psi(x)|, \|D^2\psi(x)\| \leqslant m^*, \quad \forall x \in \overline{\Omega_1}\setminus
\widetilde{\omega}_1.$$
The function $\beta$ in Proposition \ref{Prop:weight} will be of the form
\begin{equation} \label{def:beta}
\beta(x):=e^{2\lambda\|\psi\|_{L^\infty(\Omega_1)}} - e^{\lambda \psi(x)}
\end{equation}
for an appropriate parameter $\lambda>0$. From (\ref{def:beta}), we get
$$\nabla\beta(x)=-\lambda\nabla\psi(x)e^{\lambda\psi(x)},$$
$$D^2\beta(x)=-\Big(\lambda^2 \nabla \psi(x) \otimes \nabla \psi(x) + \lambda D^2\psi(x) \Big)e^{\lambda\psi(x)},$$
$$\Delta \beta(x)=-\Big(\lambda^2 |\nabla\psi(x)|^2+\lambda\Delta\psi(x) \Big)e^{\lambda\psi(x)}.$$
Using the above relations we get, for any $x \in \overline{\Omega_1}\setminus
\widetilde{\omega}_1$,
\begin{equation} \label{calcul_ineq1}
\begin{array}{ll}
  & (1-a) (\Delta \beta) |Z|^2 - 2 D^2 \beta(Z,Z) \\
= & \Big( (a-1) [\lambda^2 |\nabla\psi|^2+\lambda\Delta\psi] |Z|^2 
+ 2 [\lambda^2 \langle \nabla \psi,Z\rangle^2  + \lambda D^2\psi(Z,Z)]
\Big) e^{\lambda\psi} \\
= & \lambda^2 \Big( (a-1)|\nabla\psi|^2 |Z|^2 + 2 \langle \nabla \psi,Z\rangle^2 \Big) e^{\lambda\psi}  
+ \lambda \Big( (a-1)\Delta\psi|Z|^2 + 2 D^2\psi(Z,Z)  \Big) e^{\lambda\psi} \\
\geqslant &  \Big( \lambda^2 (a-1) m_*^2 - \lambda (a+1)m^*  \Big)  |Z|^2  
\end{array}
\end{equation}
and
\begin{equation}\label{calcul_ineq2}
\begin{array}{ll}
 & (a-1) (\Delta \beta) |\nabla \beta|^2 - 2 D^2 \beta(\nabla\beta,\nabla\beta) \\
=& \Big[   (1-a)\Big( \lambda^2 |\nabla\psi|^2+\lambda\Delta\psi \Big) \lambda^2 |\nabla\psi|^2  + 2  \Big(  \lambda^2 |\nabla\psi|^4 + \lambda D^2\psi(\nabla\psi,\nabla\psi) \Big) 
   \Big]e^{3\lambda\psi}\\
=& \Big[ \lambda^4 (3-a)|\nabla\psi|^4 + \lambda^3 \Big( (1-a)\Delta\psi|\nabla\psi|^2+2D^2\psi(\nabla\psi,\nabla\psi) \Big) \Big] e^{3\lambda\psi} \\
\geqslant & 
\Big( \lambda^4 (3-a) m_*^4 - \lambda^3 (a+1)(m^*)^3 \Big).
\end{array}
\end{equation}
The conclusion follows taking, for example,
$$\lambda:=\max\left\{ \frac{2(a+1)m^*}{(a-1)m_*^2} ; \frac{2(a+1)(m^*)^3}{(3-a)m_*^4} \right\},$$ 
$$C_1:=\frac{(a-1)m_*^2 \lambda^2}{2}, \quad C_2:=\frac{(3-a)m_*^4 \lambda^4}{2}. \hfill \Box$$

\subsection{Proof of the Carleman inequality}
\label{subsec:Carlm_proof}

Let $\tilde{\omega}_1$ be an open subset such that $\tilde{\omega}_1 \subset \omega_1$.
All the computations of the proof will be made by assuming first that
$u \in H^1(0,T;L^2(\Omega_1)) \cap L^2(0,T;H^2 \cap H^1_0(\Omega_1))$.
Then, the conclusion will follow by a density argument.
\\

Let $a \in (1,3)$ and $\beta$ be as in Proposition \ref{Prop:weight}. Let us consider the weight function
\begin{equation} \label{def-alpha}
\alpha(t,x):=\frac{\beta (x)}{t(T-t)}\, ,\quad (t,x) \in (0,T) \times \mathbb{R}^{N_1}\,,
\end{equation}
and set
\begin{equation} \label{def-z}
z(t,x):=u(t,x) e^{-M \alpha(t,x)}\, ,
\end{equation}
where $M = M(T,\mu,\beta)>0$ will be chosen later on.
One has
\begin{equation} \label{P123}
e^{-M \alpha} \mathcal{P}_{\mu,\gamma} u = P_{1}z + P_{2} z + P_{3} z\, ,
\end{equation}
where
\begin{equation} \label{def:P1-P2-P3}
\begin{array}{c}
P_{1}z:= - \Delta z +(M\alpha_{t}-M^2 |\nabla \alpha|^{2}) z + \epsilon \mu |x|^{2\gamma} b(x) z \, , \\
P_{2}z:= \frac{\partial z}{\partial t} - 2 M \nabla \alpha.\nabla z - a M (\Delta \alpha )z \, , \\
P_{3}z:= (a-1) M (\Delta \alpha) z + (1-\epsilon) \mu |x|^{2\gamma} b(x) z \, ,
\end{array}
\end{equation}
and $\epsilon \in \{0,1\}$ will be chosen later on.
We develop the classical proof (see \cite{Fursikov-Imanuvilov-186}), taking the $L^{2}(Q)$-norm in the identity (\ref{P123}),
then developing the double product, which leads to
\begin{equation} \label{P1P2<P3}
\int_{Q} \left( P_{1}z P_{2}z  - \frac{1}{2} |P_{3}z|^{2} \right) dx dt \leqslant 
\frac{1}{2} \int_{Q} | e^{-M \alpha} \mathcal{P}_{\mu,\gamma} u  |^2 dx dt\, ,
\end{equation}
where $Q:=(0,T) \times \Omega_1$ and we compute precisely each term, paying attention to the behaviour of the different
constants with respect to $\mu$ and $T$.
\\

\noindent \textbf{Step 1: Computation of the terms in (\ref{P1P2<P3}).}
\\

\noindent \textbf{Terms concerning $-\Delta z$ in $P_1 P_2$:} Integrating by parts, we get
\begin{equation} \label{Carlm1}
- \int_Q \Delta z \frac{\partial z}{\partial t} dx dt = \int_Q \nabla z\cdot
\frac{\partial}{\partial t} \nabla z dx dt =
\frac{1}{2} \int_0^T \frac{d}{dt} \int_\Omega |\nabla z|^2 dx dt =0
\end{equation}
because $\partial_t z(t,.)=0$ on $\partial \Omega_1$ and $z(0,.)=z(T,.)=0$, by definition of $z$.
Using Green's formula and the relation $\nabla z = \frac{\partial z}{\partial \nu} \nu$ on $\partial \Omega_1$, we get
\begin{equation} \label{Carlm2}
\begin{array}{ll} 
& 2 M \int_Q \Delta z \nabla \alpha \cdot \nabla z dx dt
\\ = &
2 M \int_0^T \int_{\partial \Omega_1}  \frac{\partial z}{\partial \nu} 
\nabla z \cdot \nabla \alpha d\sigma dt
- 2M \int_Q \Big( \alpha_{j,k} z_j z_k + \frac{1}{2}\nabla ( |\nabla z|^2 )
\cdot \nabla \alpha \Big) dx dt 
\\ = &
2 M \int_0^T \int_{\partial \Omega_1} \left( \frac{\partial z}{\partial \nu} \right)^2 \frac{\partial \alpha}{\partial \nu} d\sigma dt
+ M \int_Q \Big( |\nabla z|^2 \Delta \alpha - 2 D^2 \alpha \cdot(\nabla z,\nabla z) \Big) dx dt,
\end{array}
\end{equation}
with the suming index convention. Using Green's formula and $z=0$ on $\partial \Omega_1$, we get
\begin{equation} \label{Carlm3}
\begin{array}{ll}
a M \int_Q (\Delta \alpha) (\Delta z) z dx dt
& = - a M \int_Q \Big( |\nabla z|^2 \Delta \alpha + z \nabla z\cdot
 \nabla(\Delta \alpha)  \Big) dx dt\\
& = a M \int_Q \Big( - (\Delta \alpha) |\nabla z|^2  + \frac{1}{2} (\Delta^2 \alpha) z^2 \Big) dx dt.
\end{array}
\end{equation}

\noindent \textbf{Terms concerning $(M\alpha_{t}-M^2 |\nabla \alpha|^{2}) z$ in $P_1 P_2$:} Integrating by parts, we get
\begin{equation} \label{Carlm4}
\int_Q (M\alpha_{t}-M^2 |\nabla \alpha|^{2}) z \frac{\partial z}{\partial t} dx dt
= - \frac{1}{2} \int_Q (M\alpha_{t}-M^2 |\nabla \alpha|^{2})_t |z|^2 dx dt.
\end{equation}
Using the Green formula and $z(t,.)=0$ on $\partial \Omega_1$ we get
\begin{equation} \label{Carlm5}
\begin{array}{ll}
& -2M\int_Q (M\alpha_{t}-M^2 |\nabla \alpha|^{2}) z \nabla \alpha\cdot
\nabla z dx dt 
\\ = &
M \int_Q \text{div}[  (M\alpha_{t}-M^2 |\nabla \alpha|^{2}) \nabla \alpha ] |z|^2 dx dt.
\end{array}
\end{equation}
Finally, the last term concerning $(M\alpha_{t}-M^2 |\nabla \alpha|^{2}) z$ is
\begin{equation} \label{Carlm6}
-aM \int_Q (M\alpha_{t}-M^2 |\nabla \alpha|^{2}) (\Delta \alpha) |z|^2 dx dt.
\end{equation}

\noindent \textbf{Terms concerning $\epsilon \mu |x|^{2\gamma} b(x) z$ in $P_1 P_2$:} Integrating by parts, we get
\begin{equation} \label{Carlm7}
\int_Q  \epsilon \mu |x|^{2\gamma} b(x) z \frac{\partial z}{\partial t} dx dt =
\int_0^T \frac{1}{2} \frac{d}{dt} \int_{\Omega_1} \epsilon \mu |x|^{2\gamma} b(x) |z|^2 dx dt =0
\end{equation}
because $z(0,.)=z(T,.)=0$. Using Green's formula and $z(t,.)=0$ on $\partial \Omega_1$, we get
\begin{equation} \label{Carlm8}
- 2 \epsilon M \mu \int_Q  |x|^{2\gamma} b(x) z (\nabla \alpha)\cdot
(\nabla z) dx dt
= \epsilon M \mu \int_Q \text{div}[ |x|^{2\gamma} b(x) \nabla \alpha] |z|^2 dx dt.
\end{equation}
Finally, the last term concerning $\epsilon \mu |x|^{2\gamma} b(x) z$ is
\begin{equation} \label{Carlm9}
- \epsilon a M \mu \int_Q |x|^{2\gamma} b(x) (\Delta \alpha) |z|^2 dx dt.
\end{equation}

\noindent \textbf{Terms concerning $P_3$:} We have
$$\begin{array}{ll}
& -\frac{1}{2} \int_Q |P_{3}z|^{2} dx dt = -\frac{1}{2} \int_Q \Big| (a-1) M (\Delta \alpha) z + (1-\epsilon) \mu |x|^{2\gamma} b(x) z \Big|^2 dx dt
\\ \geqslant &
-\int_Q \Big( (a-1)^2 M^2 (\Delta \alpha)^2 |z|^2 + (1-\epsilon)^2 \mu^2 |x|^{4\gamma} b(x)^2 |z|^2 \Big) dx dt.
\end{array}$$

Using the relations (\ref{Carlm1}) to (\ref{Carlm9}), 
the inequality (\ref{P1P2<P3}) may be written
\begin{equation} \label{In-Carl-1}
\begin{array}{l}
\int\limits_0^T \int\limits_{\partial \Omega_1} 2M \left| \frac{\partial z}{\partial \nu} \right|^2 \frac{\partial \alpha}{\partial \nu}
+ \int\limits_Q M \Big( (1-a) (\Delta \alpha)|\nabla z|^2 - 2 D^2 \alpha .(\nabla z,\nabla z) \Big) 
\\
+ \int\limits_Q \Big\{ f + \epsilon M\mu\left( \text{div}\left[ |x|^{2\gamma} b(x) \nabla \alpha \right]  - a |x|^{2\gamma} b(x) \Delta \alpha   \right) 
-(1-\epsilon)^2 \mu^2 |x|^{4\gamma} b(x)^2 \Big\} |z|^2 
\\
\leqslant \frac{1}{2} \int\limits_Q \left| \mathcal{P}_{\mu,\gamma} u \right|^2 e^{-2M\alpha}\, ,
\end{array}
\end{equation}
where
$$\begin{array}{ll}
f:= &
\frac{Ma}{2} (\Delta^2 \alpha)
- \frac{1}{2} \left( M\alpha_t - M^2 |\nabla \alpha|^2 \right)_t
+ M \text{div} \left[ (M\alpha_t - M^2|\nabla \alpha|^2)\nabla \alpha \right]
\\ &
-aM\left( M\alpha_t - M^2 |\nabla \alpha|^2 \right) \Delta \alpha
- (a-1)^2 M^2 (\Delta\alpha)^2 
\end{array}$$

\noindent \textbf{Step 2: Estimation of the terms in (\ref{In-Carl-1}).}
\\

Using (\ref{beta_bord}), we deduce that
\begin{equation} 
\begin{array}{l}
\int\limits_Q M \Big( (1-a) (\Delta \alpha)|\nabla z|^2 - 2 D^2\alpha
.(\nabla z,\nabla z) \Big)dxdt 
\\
+ \int\limits_Q \Big\{ f + \epsilon M \mu \Big( \text{div}\left[ |x|^{2\gamma} b(x) \nabla \alpha \right]  - a |x|^{2\gamma} b(x) \Delta \alpha   \Big) 
-(1-\epsilon)^2 \mu^2 |x|^{4\gamma} b(x)^2 \Big\} |z|^2
\\ 
\leqslant \frac{1}{2} \int\limits_Q \left| \mathcal{P}_{\mu,\gamma} u \right|^2 e^{-2M\alpha}.
\end{array}
\end{equation}
We remark that
\begin{equation} \label{f_expand}
\begin{array}{ll}
f=& \frac{1}{[t(T-t)]^3} \Big\{ M^3 \Big[ (a-1)|\nabla\beta|^2\Delta\beta - 2 D^2\beta(\nabla\beta,\nabla\beta) \Big] \\
  &     + M^2 \Big[ (2t-T) \Big( |\nabla\beta|^2 + \text{div}(\beta\nabla\beta)- a \beta\Delta\beta \Big) - (a-1)^2 t(T-t) (\Delta\beta)^2  \Big] \\
  &     + M   \Big[ \frac{a}{2} \Delta^2\beta (t(T-t))^2 - \beta (T^2 - 3Tt + 3t^2)  \Big] \Big\}.
\end{array}
\end{equation}
Using (\ref{beta_13}) and the $C^4$ regularity of $\beta$ on the compact set $\overline{\Omega_1}$, 
we get constants $C_2=C_2(\beta), C_4=C_4(\beta), c=c(\beta)>0$ such that
$$f \geqslant \frac{1}{[t(T-t)]^3} \Big[ C_3 M^3 - c(T+T^2)M^2 - c(T+T^2)^2 M \Big],\, \forall (t,x) \in (0,T) \times (\overline{\Omega_1}-\tilde{\omega}_1),$$
$$|f| \leqslant \frac{1}{[t(T-t)]^3} \Big[ C_4 M^3 + c(T+T^2)M^2 + c(T+T^2)^2 M \Big],\, \forall (t,x) \in (0,T) \times \overline{\tilde{\omega}_1},$$
$$|(1-a) (\Delta \alpha)|\nabla z|^2 - 2 D^2\alpha.(\nabla z,\nabla z)| \leqslant \frac{C_2 M}{t(T-t)}|\nabla z|^2, \forall (t,x) \in [0,T] \times \overline{\tilde{\omega}_1}.$$
Thus, there exists $m_1=m_1(\beta), C_3'=C_3'(\beta), C_4'=C_4'(\beta)>0$ such that, for every $M \geqslant M_1(T,\beta)$
$$\begin{array}{l}
f \geqslant \frac{C_3' M^3}{[t(T-t)]^3}\, , \forall (t,x) \in (0,T) \times (\Omega_1\setminus\widetilde{\omega}_1),
\end{array}$$
$$\begin{array}{l}
|f| \leqslant  \frac{C_4' M^3}{[t(T-t)]^3}\, ,\forall (t,x) \in (0,T) \times \tilde{\omega}_1,
\end{array}$$
where
\begin{equation} \label{def:M1}
M_1(T,\beta):=m_1(\beta)(T+T^2).
\end{equation}
By (\ref{beta_13}) and the previous inequalities, we get, for $M \geqslant M_1(T,\beta)$,
\begin{equation} \label{In-Carl-2}
\begin{array}{l} 
\int\limits_{0}^{T} \int\limits_{\Omega_1\setminus\widetilde{\omega}_1}
\frac{C_{1} M}{t(T-t)} |\nabla z|^{2} dxdt
\\ 
+ \int\limits_{0}^{T} \int\limits_{\Omega_1\setminus\widetilde{\omega}_1}
\left[
\frac{C_{3}' M^{3}}{(t(T-t))^{3}} + 
\epsilon M \mu \Big( \text{div}\left[ |x|^{2\gamma} b \nabla \alpha \right]  - a |x|^{2\gamma} b \Delta \alpha   \Big)
-(1-\epsilon)^2 \mu^2 |x|^{4\gamma} b^2 \right] |z|^{2}
\\
\leqslant 
\int\limits_{0}^{T} \int\limits_{\tilde{\omega}_1}
\frac{C_{2} M}{t(T-t)}  |\nabla z|^2 \\
+ \int\limits_{0}^{T} \int\limits_{\tilde{\omega}_1}
\left[ 
\frac{C_{4}' M^{3}}{(t(T-t))^{3}} 
- \epsilon M\mu \Big( \text{div}\left[ |x|^{2\gamma} b \nabla \alpha \right]  - a |x|^{2\gamma} b \Delta \alpha   \Big)
+ (1-\epsilon)^2 \mu^2 |x|^{4\gamma} b^2
\right]|z|^2
\\ 

+ \frac{1}{2} \int\limits_{Q} | e^{-M \alpha} \mathcal{P}_n u  |^2 dxdt\, .
\end{array}
\end{equation}

\noindent \textbf{Step 3: End of the proof when $\gamma \in [1/2,1]$.} We take $\epsilon=1$. Then (\ref{In-Carl-2}) writes
\begin{equation} \label{In-Carl-2-gamma-gd}
\begin{array}{l} 
\int\limits_{0}^{T} \int\limits_{\Omega_1\setminus\widetilde{\omega}_1}
\frac{C_{1} M}{t(T-t)} |\nabla z|^{2} dxdt
\\ 
+ \int\limits_{0}^{T} \int\limits_{\Omega_1\setminus\widetilde{\omega}_1}
\left[
\frac{C_{3}' M^{3}}{(t(T-t))^{3}} + 
\epsilon M \mu \Big( \text{div}\left[ |x|^{2\gamma} b(x) \nabla \alpha \right]  - a |x|^{2\gamma} b(x) \Delta \alpha   \Big)
 \right] |z|^{2}
\\
\leqslant 
\int\limits_{0}^{T} \int\limits_{\tilde{\omega}_1}
\frac{C_{2} M}{t(T-t)}  |\nabla z|^2 \\
+ \int\limits_{0}^{T} \int\limits_{\tilde{\omega}_1}
\left[ 
\frac{C_{4}' M^{3}}{(t(T-t))^{3}} 
- \epsilon M\mu \Big( \text{div}\left[ |x|^{2\gamma} b(x) \nabla \alpha \right]  - a |x|^{2\gamma} b(x) \Delta \alpha   \Big)
\right]|z|^2
\\ 
+ \frac{1}{2} \int\limits_{Q} | e^{-M \alpha} \mathcal{P}_n u  |^2 dxdt\, .
\end{array}
\end{equation}
There exists $C_5=C_5(\beta)>0$ such that
\begin{equation} \label{diverging_term}
\Big|  M\mu\left( \text{div}\left[ |x|^{2\gamma} b(x) \nabla \alpha \right]  - a |x|^{2\gamma} b(x) \Delta \alpha   \right)
\Big|  \leqslant \frac{C_5 M \mu}{t(T-t)}, \forall (t,x) \in Q.
\end{equation}
Let $M_2=M_2(T,\mu,\beta)$ be defined by
\begin{equation} \label{def:M2}
M_2=M_2(T,\mu,\beta):=\sqrt{\frac{2C_5}{C_3'}} \sqrt{\mu} \left( \frac{T}{2} \right)^{2}\,.
\end{equation}
From now on, we take 
\begin{equation} \label{def:M}
M = M(T,\mu,\beta) := \mathcal{C}_2 \max \{ T+T^2 ; \sqrt{\mu} T^2 \}\, 
\end{equation}
where 
$$\mathcal{C}_2=\mathcal{C}_2(\beta):=\max \left\{ m_1 ; \sqrt{\frac{C_5}{8 C_3'}} \right\}$$
so that $M \geqslant M_1$ and $M_2$ (see (\ref{def:M1}) and (\ref{def:M2})). 
It is only at this step that the dependance of $M$ with respect to $\mu$ has to be specified.
From $M \geqslant M_2$, we deduce that
$$\Big|  M\mu\left( \text{div}\left[ |x|^{2\gamma} b(x) \nabla \alpha \right]  - a |x|^{2\gamma} b(x) \Delta \alpha   \right)
\Big|  \leqslant \frac{C_3' M^3 }{2[t(T-t)]^3}, \forall (t,x) \in Q.$$
Indeed,
$$\begin{array}{ll}
\frac{C_5 M \mu}{t(T-t)} & = \frac{C_3' M^3 }{2[t(T-t)]^3} \frac{2 \mu C_5 [t(T-t)]^2}{C_3' M^2} \\
                         & \leqslant \frac{C_3' M^3 }{2[t(T-t)]^3} \frac{\mu C_5 T^4}{2 C_3' M_2^2} \\
                         & \leqslant \frac{C_3' M^3 }{2[t(T-t)]^3}.
\end{array}$$
Thus, (\ref{In-Carl-2-gamma-gd}) implies
\begin{equation} \label{In-Carl-3}
\begin{array}{ll} 
\int\limits_{0}^{T} \int\limits_{\Omega_1\setminus\widetilde{\omega}_1}
\left(\frac{C_{1} M}{t(T-t)} |\nabla z|^{2} 
+ \frac{C_{3}' M^3}{2(t(T-t))^{3}} |z|^{2}\right) dxdt
\leqslant 
& 
\int\limits_{0}^{T} \int\limits_{\tilde{\omega}_1}
\left( \frac{C_{2} M}{t(T-t)}  |\nabla z|^2 
+ \frac{C_6 M^{3}}{(t(T-t))^{3}} |z|^2\right) dxdt
\\ &
+  \frac{1}{2} \int\limits_{Q} | e^{-M \alpha} \mathcal{P}_{\mu,\gamma} u  |^2 \, ,
\end{array}
\end{equation}
where $C_6=C_6(\beta):=C_4'+C_3'/2$. 
For every $\epsilon'>0$, we have
\begin{multline}\label{estInCarl}
\frac{C_{1} M}{t(T-t)} | \nabla u  - M (\nabla \alpha) u  |^{2}
+
\frac{C_{3}' M^{3}}{2 (t(T-t))^{3}} |u|^{2}
\\
\geqslant
\left( 1-\frac{1}{1+\epsilon'} \right) \frac{C_{1} M}{t(T-t)} | \nabla u |^{2}
+ \frac{M^3}{(t(T-t))^3}\left( \frac{C_{3}'}{2} - \epsilon' C_1 |\nabla \beta|^2 \right) |u|^2\,.
\end{multline}
Hence, choosing
$$\epsilon'=\epsilon'(\beta):=\frac{C_3'}{4C_1\|\nabla \beta\|_\infty^2}\, ,$$ 
from \eqref{In-Carl-3}, \eqref{estInCarl} and \eqref{def-z} we deduce that
\begin{multline} 
\int_{0}^{T} \int\limits_{\Omega_1\setminus\widetilde{\omega}_1}  
\left(
\frac{C_{7} M}{t(T-t)} |\nabla u|^{2}
+ \frac{C_{3}' M^{3}|u|^{2}}{4(t(T-t))^{3}} 
\right)  e^{-2M\alpha} dxdt
\\
\leqslant 
\int_{0}^{T} \int_{\tilde{\omega}_1} 
\left(
\frac{C_{8} M}{t(T-t)}  | \nabla u |^{2} 
+
\frac{C_{9} M^{3} |u|^{2}}{(t(T-t))^{3}}
\right)e^{-2M\alpha} 
\\+ 
\frac{1}{2} \int_{Q} | e^{-M \alpha} \mathcal{P}_{\mu,\gamma} u  |^2 \, ,
\end{multline}
where 
$C_7=C_7(\beta):=[1-1/(1+\epsilon')]C_1$, 
$C_{8}=C_{8}(\beta):=2C_{2}$ and 
$C_{9}=C_{9}(\beta):=C_{6}+2C_{2} \|\nabla\beta\|_\infty^2$.
So, adding the same quantity to both sides, 
\begin{multline} \label{In-Carl-4}
\int_Q  
\left(
\frac{C_{7} M}{t(T-t)} |\nabla u|^2 
+ \frac{C_{3}' M^{3}|u|^{2}}{4(t(T-t))^{3}} 
\right)  e^{-2M \alpha}  \leqslant 
\frac{1}{2} \int_{Q} | e^{-M \alpha} \mathcal{P}_{\mu,\gamma} u  |^2 
\\
+
\int_{0}^{T} \int_{\tilde{\omega}_1} 
\left(
\frac{C_{10} M}{t(T-t)}  |\nabla u|^2 
+
\frac{C_{11} M^{3} |u|^{2}}{(t(T-t))^{3}} 
\right)e^{-2M\alpha} \, ,
\end{multline}
where $C_{10}=C_{10}(\beta):=C_{8}+C_7$ and $C_{11}=C_{11}(\beta):=C_{9}+C_3'/4$.
Let us prove that the third term on the right-hand side can be dominated by 
terms similarly to the other two ones.
We consider $\rho \in C^{\infty}(\mathbb{R}^{N_1};\mathbb{R}_{+})$ such that $0\le \rho \le 1$ and
$$
\rho \equiv 1 \text{  on  } \tilde{\omega}_1\ \text{  and  }
\rho \equiv 0 \text{  on  } \mathbb{R}^{N_1}-\omega_1.
$$
We have
$$\int_{Q} (\mathcal{P}_{\mu,\gamma} u) \frac{u \rho e^{-2M\alpha}}{t(T-t)}
dxdt 
= \int_0^T \int_{\Omega_1} \left[  
\frac{\partial u}{\partial t} - \Delta u  + \mu |x|^{2\gamma} b(x) u
\right]
\frac{u \rho e^{-2M\alpha}}{t(T-t)}.$$
Integrating by parts with respect to time and space, we obtain
$$\int_Q \frac{1}{2} \frac{\partial (u^2)}{\partial t} 
\frac{\rho e^{-2M\alpha}}{t(T-t)} dxdt
=\int_Q \frac{1}{2} |u|^2 \rho \left(
\frac{2M\alpha_t}{t(T-t)} + \frac{T-2t}{(t(T-t))^2}
\right) e^{-2M\alpha} $$
and
\begin{multline}
-\int_Q \Delta u \frac{u \rho e^{-2M\alpha}}{t(T-t)} dxdt
= \int_Q 
\frac{\rho e^{-2M\alpha}}{t(T-t)} |\nabla u|^2  
\\ 
- \int_Q  \frac{|u|^2 e^{-2M\alpha}}{2t(T-t)} 
\left( \Delta \rho - 4 M \nabla \rho.\nabla \alpha + \rho (4M^2|\nabla\alpha|^2 -2M\Delta\alpha) \right)\,.
\end{multline}
Thus,
\begin{multline}
\int_{Q} \mathcal{P}_{\mu,\gamma} u \frac{u \rho e^{-2M\alpha}}{t(T-t)} dxdt
\geqslant
\int_Q \frac{\rho e^{-2M\alpha}}{t(T-t)} |\nabla u |^2   
\\
- \int_Q  \frac{|u|^2 e^{-2M\alpha}}{2t(T-t)} 
\left( \Delta \rho - 4 M \nabla \rho\cdot\nabla \alpha + \rho \left( 4M^2 |\nabla \alpha|^2 - 2M\Delta \alpha - 2M\alpha_t-\frac{T-2t}{t(T-t)} \right) \right) \,.
\end{multline}
Therefore,
\begin{multline*}
\int_0^T \int_{\tilde{\omega}_1} \frac{C_{10} M}{t(T-t)} |\nabla u|^2   
e^{-2M\alpha} dxdt
\\ \leqslant
\int_Q \frac{C_{10} M \rho}{t(T-t)} |\nabla u|^2  e^{-2M\alpha} 
 \leqslant
\int_{Q}
 \mathcal{P}_{\mu,\gamma} u \frac{C_{10} M u \rho e^{-2M\alpha}}{t(T-t)} 
\\
+ \int_{Q}
\frac{C_{10} M |u|^2 e^{-2M\alpha}}{2t(T-t)} 
\left( \Delta \rho - 4 M \nabla \rho\cdot\nabla \alpha + \rho \left( 4M^2|\nabla\alpha|^2-2M\Delta \alpha -2M\alpha_t-\frac{T-2t}{t(T-t)} \right) \right)
\\ \leqslant
\int_Q |\mathcal{P}_{\mu,\gamma} u|^2 e^{-2M\alpha} +
\int_0^T \int_{\omega_1} \frac{C_{12} M^3 |u|^2 e^{-2M\alpha}}{(t(T-t))^3}
dxdt 
\end{multline*}
for some constant $C_{12}=C_{12}(\beta,\rho)>0$. Combining (\ref{In-Carl-4}) with the previous inequality, we get
\begin{multline} 
\int_Q  
\left(
\frac{C_{7} M}{t(T-t)} |\nabla u|^2 
+ \frac{C_{3}' M^{3}|u|^{2}}{4(t(T-t))^{3}} 
\right)  e^{-2M\alpha} dxdt 
\\
\leqslant 
 \int_{Q} 2 | e^{-M\alpha} \mathcal{P}_{\mu,\gamma} u  |^2 
+
\int_{0}^{T} \int_{\omega_1} \frac{C_{13} M^{3} |u|^{2}}{(t(T-t))^{3}} e^{-2M\alpha}\, ,
\end{multline}
where $C_{13}=C_{13}(\beta,\rho):=C_{11}+C_{12}$. 
Then, the global Carleman estimates (\ref{Carl_est}) holds with
$$\mathcal{C}_1=\mathcal{C}_1(\beta):=\frac{\min\{C_7  ; C_3' /4 \}}{\max\{ 2 ; C_{13}  \}}.$$

\noindent \textbf{Step 4: End of the proof when $\gamma \in (0,1/2)$.}  
The left-hand side of (\ref{diverging_term}) diverges at $x=0$, thus the proof cannot be ended in the same way and
we take $\epsilon=0$. Then (\ref{In-Carl-2}) writes
\begin{equation} \label{In-Carl-2-gamma-pt}
\begin{array}{l} 
\int\limits_{0}^{T} \int\limits_{\Omega_1\setminus\widetilde{\omega}_1}
\frac{C_{1} M}{t(T-t)} |\nabla z|^{2} dxdt
\\ 
+ \int\limits_{0}^{T} \int\limits_{\Omega_1\setminus\widetilde{\omega}_1}
\left[
\frac{C_{3}' M^{3}}{(t(T-t))^{3}} -(1-\epsilon)^2 \mu^2 |x|^{4\gamma} b(x)^2 \right] |z|^{2} dxdt
\\
\leqslant 
\int\limits_{0}^{T} \int\limits_{\tilde{\omega}_1}
\frac{C_{2} M}{t(T-t)}  |\nabla z|^2 \\
+ \int\limits_{0}^{T} \int\limits_{\tilde{\omega}_1}
\left[ 
\frac{C_{4}' M^{3}}{(t(T-t))^{3}} 
+ (1-\epsilon)^2 \mu^2 |x|^{4\gamma} b(x)^2
\right]|z|^2
\\ 

+ \frac{1}{2} \int\limits_{Q} | e^{-M \alpha} \mathcal{P}_n u  |^2 dxdt\, .
\end{array}
\end{equation}
Let
\begin{equation} \label{def:M2_bis}
M_2=M_2(T,\beta,\mu):=\frac{T^2}{4} \mu^{2/3} \sqrt[3]{\frac{2R^{4\gamma}\|b\|_{\infty}^2}{C_3'}}
\end{equation}
where $R>0$ is such that $\Omega_1 \subset B(0,R)$. 
From now on, we take 
\begin{equation} \label{def:M}
M = M(T,\mu,\beta) := \mathcal{C}_2 \max \{ T+T^2 ; \mu^{2/3} T^2 \}\, 
\end{equation}
where 
$$\mathcal{C}_2=\mathcal{C}_2(\beta):=\max \left\{ m_1 ; \frac{1}{4} \sqrt[3]{\frac{2R^{4\gamma}\|b\|_\infty^2}{C_3'}} \right\}$$
so that $M \geqslant M_1$ and $M_2$ (see (\ref{def:M1}) and (\ref{def:M2_bis})). 
It is only at this step that the dependence of $M$ with respect to $\mu$ has to be specified.
From $M \geqslant M_2$, we deduce that
$$\mu^2 |x|^{4\gamma} b(x)^2  \leqslant \frac{C_{3}' M^{3}}{2(t(T-t))^{3}}, \forall (t,x) \in Q.$$
From (\ref{In-Carl-2-gamma-pt}), we are lead to an inequality of the form (\ref{In-Carl-3}) and the proof may be finished as 
in Step 3. $\hfill \Box$

\bibliography{biblio3_BIS}

\begin{thebibliography}{10}

\bibitem{Grushin}
K.~Beauchard, P.~Cannarsa, and R.~Guglielmi.
\newblock {Null controllability of Grushin-type equations in dimension two}.
\newblock {\em JEMS (to appear)}, 2012.

\bibitem{BK}
A.~L. Bukhgeim and M.~V. Klibanov.
\newblock Global uniqueness of a class of multidimensional inverse problems.
\newblock {\em Sov. Math. Dokl.}, 24:244--247, 1981.

\bibitem{CTY1}
P.~Cannarsa, J.~Tort, and M.~Yamamoto.
\newblock Determination of source terms in a degenerate parabolic equation.
\newblock {\em Inverse Problems}, 26:105003, 2010.

\bibitem{CTY2}
P.~Cannarsa, J.~Tort, and M.~Yamamoto.
\newblock Unique continuation and approximate controllability for a degenerate
  parabolic equation.
\newblock {\em Applicable Analysis}, 91:1409--1425, 2012.

\bibitem{capogna}
Luca Capogna, Donatella Danielli, Scott~D. Pauls, and Jeremy~T. Tyson.
\newblock {\em An introduction to the Heisenberg group and the sub-Riemannian
  isoperimetric problem}, volume Progress in Mathematics 259.
\newblock Birkhauser Verlag, Basel, 2007.

\bibitem{JMC-book}
J.-M. Coron.
\newblock {\em Control and nonlinearity}, volume 136.
\newblock Mathematical Surveys and Monographs, 2007.

\bibitem{Fursikov-Imanuvilov-186}
A.V. Fursikov and O.Y. Imanuvilov.
\newblock Controllability of evolution equations.
\newblock {\em Lecture Notes Series, Seoul National University Research
  Institute of Mathematics Global Analysis Research Center, Seoul}, 34, 1996.

\bibitem{IY1}
O.~Y. Imanuvilov and M.~Yamamoto.
\newblock {Lipschitz stability in inverse parabolic problems by the Carleman
  estimate}.
\newblock {\em Inverse Problems}, 14:1229--45, 1998.

\bibitem{Is}
V.~Isakov.
\newblock Inverse source problems.
\newblock 1990.

\bibitem{Kl}
M.~V. Klibanov.
\newblock {Inverse problems and Carleman estimates}.
\newblock {\em Inverse Problems}, 8:575--96, 1992.

\bibitem{Lebeau-Robbiano}
G.~Lebeau and L.~Robbiano.
\newblock Contr{\^o}le exact de l'{\'e}quation de la chaleur.
\newblock {\em Comm. P.D.E.}, 20:335--356, 1995.

\bibitem{Lebeau-LeRousseau}
G.~Lebeau and J.~Le Rousseau.
\newblock {On Carleman estimates for elliptic and parabolic operators.
  Applications to unique continuation and control of parabolic equations}.
\newblock {\em ESAIM:COCV (DOI:10.1051 cocv 2011168)}, 2011.

\bibitem{resi}
M.~Reed and B.~Simon.
\newblock {\em Methods of modern mathematical physics. I. Functional analysis.
  Second edition.}
\newblock Academic Press, Inc. [Harcourt Brace Jovanovich, Publishers], 1980.

\bibitem{Weyl}
H.~Weyl.
\newblock {\em The theory of groups and quantum mechanics}.
\newblock Methuen, 1931.

\bibitem{Ya}
M.~Yamamoto.
\newblock Carleman estimates for parabolic equations and applications.
\newblock {\em Inverse Problems}, 25:123013, 2009.

\end{thebibliography}
\bibliographystyle{plain}
\end{document}